\documentclass[12pt,draft,leqno]{article}
\usepackage{amssymb, eucal, latexsym}

\newtheorem{theorem}{Theorem}[section]
\newtheorem{lm}[theorem]{Lemma}
\newtheorem{exa}[theorem]{Example}

\newtheorem{cor}[theorem]{Corollary}
\newtheorem{pro}[theorem]{Proposition}
\newtheorem{defi}[theorem]{Definition}
\newtheorem{defis}[theorem]{Definitions}

\newtheorem{notas}[theorem]{Notations}
\newtheorem{rem}[theorem]{Remark}
\newtheorem{rems}[theorem]{Remarks}
\newtheorem{fact}[theorem]{Fact}

\newtheorem{nist}[theorem]{}

\def\p{\varphi}
\def\a{\alpha}

\def\d{\delta}

\def\g{\gamma}
\def\GA{\Gamma}

\def\l{\lambda}
\def\LAM{\Lambda}

\def\TE{\Theta}

\def\s{\sigma}
\def\SI{\Sigma}

\def\pl{\varphi_\Lambda}

\def\lag{\lambda_A^g}
\def\lbg{\lambda_B^g}

\def\lra{\longrightarrow}

\def\dar{\downarrow}

\def\sbe{\subseteq}
\def\spe{\supseteq}
\def\stm{\setminus}
\def\ems{\emptyset}
\def\nes{\neq\emptyset}

\def\cuk{\,\check{}\,}

\def\fa{\forall}

\def\we{\wedge}
\def\bw{\bigwedge}
\def\bv{\bigvee}

\def\ap{^\prime}
\def\inv{^{-1}}
\def\st{\ |\ }

\def\llx{\ll_{\rho}}
\def\lle{\ll_{\eta}}

\def\nin{\not\in}

\def\card #1{\vert #1 \vert}

\def\OO{{\cal O}}

\def\SKLC{{\bf SkeLC}}
\def\SKAL{{\bf DSkeLC}}

\def\SAC{{\bf DSkeC}}

\def\OLC{{\bf OpLC}}
\def\OAL{{\bf DOpLC}}

\def\HC{{\bf HC}}

\def\PLC{{\bf PLC}}

\def\SkeZLC{{\bf SkeZLC}}
\def\SkeZLBA{{\bf SkeZLBA}}

\def\HLC{{\bf HLC}}
\def\GBPL{{\bf GBPL}}
\def\ZLC{{\bf ZLC}}
\def\DInZLC{{\bf DInZLC}}
\def\InZLC{{\bf InZLC}}
\def\DLC{{\bf DLC}}
\def\PDLC{{\bf PDLC}}

\def\QZHC{{\bf QZHC}}

\def\Bool{{\bf Bool}}
\def\BoolC{{\bf BoolC}}
\def\OBool{{\bf OBool}}
\def\ZHC{{\bf ZHC}}
\def\OZHC{{\bf OZHC}}
\def\PZLC{{\bf PZLC}}
\def\OPZLC{{\bf OPZLC}}
\def\OZLC{{\bf OZLC}}
\def\QPZLC{{\bf QPZLC}}
\def\QPZLBA{{\bf QPZLBA}}

\def\PZLBA{{\bf PZLBA}}
\def\ZLBA{{\bf ZLBA}}
\def\OZLBA{{\bf OZLBA}}
\def\OPZLBA{{\bf OPZLBA}}
\def\QGBPL{{\bf QGBPL}}
\def\OGBPL{{\bf OGBPL}}

\def\OLC{{\bf OpLC}}
\def\OAL{{\bf DOpLC}}
\def\SAC{{\bf DSkeC}}

\def\SKAL{{\bf DSkeLC}}
\def\SKLC{{\bf SkeLC}}

\def\B{\mbox{{\boldmath $B$}}}
\def\C{\mbox{{\boldmath $C$}}}

\def\2{\mbox{{\bf 2}}}
\def\3{\mbox{{\bf 3}}}

\def\int{\mbox{{\rm int}}}
\def\Fr{\mbox{{\rm Fr}}}
\def\cl{\mbox{{\rm cl}}}

\def\doc{\hspace{-1cm}{\em Proof.}~~}
\def\sq{\hspace*{\fill} \hbox{\vrule\vbox{\hrule\phantom{o}\hrule}\vrule}}
\def\sqs{\sq \vspace{2mm}}

\def\BBBB{{\rm I}\!{\rm B}}

\title{{\LARGE\bf
A De Vries-type Duality Theorem for Locally Compact Spaces -- III}\\
\vspace{0.35cm}
{\large\bf Georgi Dimov}\thanks{This paper was supported by the
project no. 005/2009 $``$General and Categorical Topology" of the Sofia
University $``$St. Kl. Ohridski".}\\
\vspace{0.25cm}
 {\footnotesize Dept. of Math. and
Informatics, Sofia University,  Blvd. J. Bourchier 5, 1164 Sofia,
Bulgaria}}

\author{}

\date{}

\begin{document}
\maketitle
\begin{abstract}
{\footnotesize
\noindent This paper is a continuation of the papers \cite{Di3,Di4} and also of \cite{Di5,Di2}.
In it we prove some new Stone-type duality theorems for some subcategories of the category $\ZLC$ of
locally compact zero-dimensional Hausdorff spaces and continuous maps. These theorems are new even in the compact case.
They concern the cofull subcategories $\SkeZLC$, $\QPZLC$, $\OZLC$ and $\OPZLC$ of the category $\ZLC$ determined,
respectively, by the skeletal maps, by the quasi-open perfect maps, by the open maps and by the open perfect maps. In
this way, the  zero-dimensional analogues of Fedorchuk Duality Theorem \cite{F2} and its generalization (presented in \cite{Di2}) are obtained.
Further, we characterize the injective and surjective   morphisms of the category  $\HLC$ of locally compact Hausdorff spaces and continuous maps, as well as of the category $\ZLC$,  and of their subcategories  discussed in \cite{Di2, Di3, Di4}, by means of some properties of their dual morphisms. This generalizes some well-known results of M. Stone  \cite{ST} and de Vries  \cite{dV2}. An analogous problem is investigated for the homeomorphic embeddings, dense embeddings, LCA-embeddings etc., and a generalization of a theorem of  Fedorchuk  \cite[Theorem 6]{F2} is obtained. Finally, in analogue to some well-known results of M. Stone \cite{ST}, the dual objects of the open, regular open, clopen, closed, regular closed etc. subsets of a space $X\in\card{\HLC}$ or $X\in\card{\ZLC}$ are described by means of the dual objects of $X$; some of these results (e.g., for regular closed sets)  are new even in the compact case.}
\end{abstract}

{\footnotesize {\em  MSC:} primary 54D45, 18A40; secondary 54C10,
 54E05.

}

\footnotetext[1]{{\footnotesize {\em E-mail address:}
gdimov@fmi.uni-sofia.bg}}

\baselineskip = \normalbaselineskip

\section*{Introduction}

This paper is a continuation of the papers \cite{Di3,Di4} and also of \cite{Di5,Di2}.
In the paper \cite{Di3} we extended de Vries Duality Theorem \cite{dV2} to the category $\HLC$ of locally compact Hausdorff spaces and continuous maps and denoted by $\DLC$ the obtained dual category and by $\LAM^t:\HLC\lra\DLC$ and $\LAM^a:\DLC\lra\HLC$ the corresponding duality functors. In the paper \cite{Di4} we extended in two ways the Stone Duality Theorem \cite{ST}: first,  to the category $\ZLC$ of locally compact zero-dimensional Hausdorff spaces and continuous maps, and, second, to the cofull subcategory $\PZLC$ of the category $\ZLC$ determined by the perfect maps (this means that the category $\PZLC$ has the same objects as the category $\ZLC$ but the morphisms of the category $\PZLC$ are the perfect maps). We denoted by $\ZLBA$ the obtained dual category of the category $\ZLC$, by $\TE^t:\ZLC\lra\ZLBA$ and $\TE^a:\ZLBA\lra\ZLC$ the corresponding duality functors, and by $\GBPL$ the obtained dual category of the category $\PZLC$, and by $\TE^t_g:\PZLC\lra\GBPL$ and $\TE^a_g:\GBPL\lra\PZLC$ the corresponding duality functors. The objects of the category $\GBPL$ are the generalized Boolean pseudolattices and the morphisms of this category are the {0}-pseudolattice homomorphisms satisfying an additional condition.

The first section of the paper contains some preliminary results.

In the second section of this paper, we prove some new Stone-type duality theorems for some subcategories of the category $\ZLC$. These theorems are new even in the compact case (see Theorems \ref{skenuldimth}, \ref{corperqo}(b),(c), \ref{thqgbpl},\ref{thzop}(b),\ref{thzopp}(b),\ref{thzoppg}). They concern the cofull subcategories $\SkeZLC$, $\QPZLC$, $\OZLC$ and $\OPZLC$ of the category $\ZLC$ determined,
respectively, by the skeletal maps (defined by
Mioduszewski and Rudolf in \cite{MR2}), by the quasi-open (defined by S. Marde\v{s}ic and P. Papic in \cite{MP}) perfect maps, by the open maps, and by the open perfect maps. Since the categories $\QPZLC$ and $\OPZLC$ are cofull subcategories simultaneously of the categories $\ZLC$ and $\PZLC$, 
we  find their images by the both functors $\TE^t$ and $\TE^t_g$
(see Corollary \ref{corperqo}(b), Theorem \ref{thqgbpl} and Corollary \ref{thzoppg}). For the compact case, these theorems give the following results:

\smallskip

(a) The category $\QZHC$ of compact zero-dimensional Hausdorff spaces and quasi-open maps is dually equivalent to the category $\BoolC$ of Boolean algebras and complete Boolean homomorphisms (see Corollary \ref{corperqo}(c)), and

\smallskip

(b) The category $\OZHC$ of compact zero-dimensional Hausdorff spaces and open maps is dually equivalent to the category $\OBool$ of Boolean algebras and Boolean homomorphisms $\p$ having lower adjoint $\psi$ (i.e. the pair $(\psi,\p)$ forms a Galois connection (see Corollary \ref{thzopps}(b)).

Let us mention also the following result (see Theorem \ref{thqgbpl}): The category $\QPZLC$ is dually equivalent to the cofull subcategory $\QGBPL$ of the category $\GBPL$ whose morphisms, in addition, preserve all meets that happen to exist.

From the mentioned above Corollary \ref{corperqo}(c) and Fedorchuk Duality Theorem \cite{F2}, we obtain the following assertion which is a special case of a much more general theorem of Monk \cite{Monk}: a Boolean homomorphism can be extended to a complete homomorphism between the corresponding minimal completions iff it is a complete homomorphism.

Note also that Theorem \ref{skenuldimth} and Corollary \ref{corperqo}(b),(c) are zero-dimensional analogues of Fedorchuk Duality Theorem \cite{F2} and its generalization presented in \cite{Di2}.

In the third section, we characterize the injective and surjective   morphisms of the categories $\ZLC$ and $\HLC$  and their subcategories  discussed in \cite{Di2, Di3, Di4,Di5} by means of some  properties of their dual morphisms. 
  If $f\in\HLC(X,Y)$,  then we use its dual map $\LAM^t(f)$, and  if $f\in \ZLC$ (resp., $f\in\PZLC$), we  use its dual $\TE^t(f)$ (resp., $\TE^t_g(f)$). When $f$ has some extra properties, e.g. if $f$ is open or quasi-open, the conditions on its dual map can be simplified and we do this. Such investigations were done by M. Stone in \cite{ST} for mappings $f\in\ZLC(X,Y)$ which are either closed embeddings or surjective maps, and by de Vries \cite{dV2} using his duality theorem. Here we generalize these results.
An analogous problem is investigated for the homeomorphic embeddings, dense embeddings, LCA-embeddings, closed embeddings.  Our Theorem \ref{lcaemb}, in which we characterize LCA-embeddings,  generalizes a theorem of  Fedorchuk  \cite[Theorem 6]{F2}.

Finally, in the fourth section, in analogue to some well-known results of M. Stone \cite{ST}, the dual objects of the open, regular open, clopen, closed, regular closed etc. subsets of a space $X\in\card{\HLC}$ or $X\in\card{\ZLC}$ are described by means of the dual objects $\LAM^t(X)$, $\TE^t(X)$ or $\TE^t_g(X)$ of $X$. Some of these results (e.g., for regular closed sets)  are new even in the compact case. For example, we prove the following (see Theorem \ref{clregclz}): Let $I$ and $J$ be generalized Boolean pseudolattices, $X=\TE^a_g(I)$ and $F=\TE^a_g(J)$; then
$F$ is homeomorphic to a regular closed subset of $X$ iff there exists a surjective 0-pseudolattice homomorphism $\p:I\lra J$ which preserves all meets that happen to exist in $I$. Hence, if $A$ and $B$ are Boolean algebras, $X=S^a(A)$ and $F=S^a(B)$ then
$F$ is homeomorphic to a regular closed subset of $X$ iff there exists a  complete epimorphism $\p:A\lra B$ (here $S^a:\Bool\lra \ZHC$ is the Stone duality functor).

In this paper we use the notations introduced in
 \cite{Di2,Di3,Di4,Di5}, as well as the notions and the results from
these papers. In the first section we recall many of them. Let us mention some more notations. If $(A,\le)$ is a poset
and $a\in A$ then  $\downarrow a$ is the set $\{b\in A\st b\le
a\}$. If $f:X\lra Y$ is a function and $M\sbe X$ then
$f_{\upharpoonright M}$ is the restriction of $f$ having $M$ as a
domain and $f(M)$ as a codomain.

\section{Preliminaries}

We will say that a subcategory $\B$\/ of a category $\C$\/ is a
{\em cofull subcategory} if $\card{\B}=\card{\C}$.

\begin{defis}\label{dvfi}{\rm (\cite{Di3})}
\rm 
Let  $\HLC$ be the category of all locally compact Hausdorff
spaces and all continuous maps between them.

Let $\DLC$ be the category whose objects are all complete LCAs and
whose morphisms are all functions $\p:(A,\rho,\BBBB)\lra
(B,\eta,\BBBB\ap)$ between the objects of $\DLC$ satisfying
conditions

\smallskip

\noindent(DLC1) $\p(0)=0$;\\
(DLC2) $\p(a\we b)=\p(a)\we \p(b)$, for all $a,b\in A$;\\
(DLC3) If $a\in\BBBB, b\in A$ and $a\llx b$, then $(\p(a^*))^*\lle
\p(b)$;\\
(DLC4) For every $b\in\BBBB\ap$ there exists $a\in\BBBB$ such that
$b\le\p(a)$;\\
\noindent(DLC5) $\p(a)=\bigvee\{\p(b)\st b\in\BBBB, b\llx a\}$,
for every $a\in A$;

\medskip

{\noindent}let the composition $``\diamond$" of two morphisms
$\p_1:(A_1,\rho_1,\BBBB_1)\lra (A_2,\rho_2,\BBBB_2)$ and
$\p_2:(A_2,\rho_2,\BBBB_2)\lra (A_3,\rho_3,\BBBB_3)$ of\/ $\DLC$
be defined by the formula
\begin{equation}\label{diamcon}
\p_2\diamond\p_1 = (\p_2\circ\p_1)\cuk,
\end{equation}
 where, for every
function $\psi:(A,\rho,\BBBB)\lra (B,\eta,\BBBB\ap)$ between two
objects of\/ $\DLC$, $\psi\cuk:(A,\rho,\BBBB)\lra
(B,\eta,\BBBB\ap)$ is defined as follows:
\begin{equation}\label{cukfcon}
\psi\cuk(a)=\bigvee\{\psi(b)\st b\in \BBBB, b\llx a\},
\end{equation}
for every $a\in A$.
\end{defis}

Note that in \cite{Di3} it was shown that condition (DLC3) in the  definition of the category $\DLC$ can be replaced with the following one:

\smallskip

\noindent(DLC3S) If $a, b\in A$ and $a\llx b$, then $(\p(a^*))^*\lle
\p(b)$.

In \cite{Di3}, the following duality theorem was proved:

\begin{theorem}\label{lccont}{\rm (\cite{Di3})}
The categories $\HLC$ and\/ $\DLC$ are dually equivalent. In more
details, let
$$\LAM^t:\HLC\lra\DLC \mbox{ and }\LAM^a:\DLC\lra\HLC$$
be the contravariant functors extending, respectively,  the Roeper
correspondences $\Psi^t:\card{\HLC}\lra\card{\DLC}$  and\/
$\Psi^a:\card{\DLC}\lra\card{\HLC}$  (see \cite{Di2} for $\Psi^t$ and $\Psi^a$) to the
 morphisms of the categories $\HLC$ and\/ $\DLC$ in
the following way:
$$\LAM^t(f)(G)=\cl(f\inv(\int(G))),$$
for every $f\in\HLC(X,Y)$ and every $G\in RC(Y)$, and
$$\LAM^a(\p)(\s\ap)\cap\BBBB=\{a\in \BBBB\st \mbox{if } b\in \BBBB
\mbox{ and } a\llx b \mbox{ then }\p(b)\in\s\ap\},$$
for every $\p\in\DLC((A,\rho,\BBBB),(B,\eta,\BBBB\ap))$ and for
every $\s\ap\in\LAM^a(B,\eta,\BBBB\ap)$; then
$$\l^g: Id_{\,\DLC}\lra\LAM^t\circ\LAM^a, \mbox{ where }
\l^g(A,\rho,\BBBB)=\l_A^g$$
for every $(A,\rho,\BBBB)\in\card\DLC$ (see (\cite{Di2}) for
the notation $\l_A^g$), and
$$t^l:Id_{\,\HLC}\lra\LAM^a\circ\LAM^t, \mbox{ where }t^l(X)=t_X$$
for every $X\in\card\HLC$ (see (\cite{Di2}) for the notation
$t_X$), are natural isomorphisms.
\end{theorem}

\begin{defis}\label{defplc}{\rm (\cite{Di5, Di3})}
\rm Let $\PLC$ be the cofull subcategory of the category $\HLC$
 whose morphisms are perfect maps.

Let $\PDLC$ be the cofull subcategory of the category $\DLC$
 whose morphisms are all $\DLC$-morphisms
$\p:(A,\rho,\BBBB)\lra
(B,\eta,\BBBB\ap)$  satisfying the following condition:

\smallskip

\noindent(PAL5) If $a\in\BBBB$ then $\p(a)\in\BBBB\ap$.
\end{defis}

\begin{theorem}\label{lcper}{\rm (\cite{Di5, Di3})}
The categories\/ $\PLC$ and\/ $\PDLC$ are dually equivalent.
More precisely, the restrictions
$\LAM^t_p:\PLC\lra\PDLC$  and $\LAM^a_p:\PDLC\lra\PLC$
of the contravariant functors $\LAM^t$ and $\LAM^a$ defined in Theorem \ref{lccont} are the
desired duality functors.
\end{theorem}

In \cite{Di3}, the following result was proved:

\begin{pro}\label{mainfed}
The categories  $\SKLC$ and\/ $\SKAL$ (see \cite{Di2} for these notations) are (non full) subcategories
of, respectively, $\HLC$ and $\DLC$. The restriction of the
contravariant functor $\LAM^a$ (respectively, $\LAM^t$) to  the
subcategory $\SKAL$ (resp., $\SKLC$) coincides with the
contravariant functor $\Psi_1^a$ (resp.,  $\Psi_1^t$) (see \cite{Di2} for these notations).
\end{pro}

\begin{defi}\label{deflba}{\rm (\cite{Di4})}
\rm A pair $(A,I)$, where $A$ is a Boolean algebra and $I$ is an
ideal of $A$ (possibly non proper) which is dense in $A$ (shortly,
dense ideal), is called a {\em local Boolean algebra} (abbreviated
as LBA).  Two LBAs $(A,I)$ and $(B,J)$ are said to be {\em isomorphic}
if there exists a Boolean isomorphism $\p:A\lra B$ such that
$\p(I)=J$.
\end{defi}

\begin{exa}\label{extrcr}
\rm Let $B$ be a Boolean algebra. Then there exist  a
smallest contact relations on $B$
defined by   $a\rho_s b$ iff $a\wedge b\neq 0$.
Note that, for $a,b\in B$, $a\ll_{\rho_s} b$ iff $a\le b$; hence
$a\ll_{\rho_s} a$, for any $a\in B$. Thus $(B,\rho_s)$ is a normal
contact algebra.
\end{exa}

\begin{lm}\label{rhoslemma}{\rm (\cite{Di4})}
If $(A,\rho_s,\BBBB)$ is an LCA  then $(A,\BBBB)$ is an LBA.
Conversely, for any LBA $(A,I)$, the triple $(A,\rho_s,I)$ is an
LCA.
\end{lm}

\begin{nist}\label{sii}
\rm
Since we follow Johnstone's terminology from \cite{J2}, we will use
the term {\em pseudolattice} for a poset having all finite
non-empty meets and joins; the pseudolattices with a bottom will
be called {\em $\{0\}$-pseudolattices}. Recall that a distributive
$\{0\}$-pseudolattice $A$ is called a {\em generalized Boolean
pseudolattice} if it satisfies the following condition:

\smallskip

\noindent(GBPL) for every $a,b,c\in A$ such that
$b\le a\le c$ there exists $x\in A$ with $a\we x=b$ and $a\vee
x=c$ (i.e., $x$ is the {\em relative complement of $a$ in the
interval} $[b,c]$).

\smallskip

Let $A$ be a distributive $\{0\}$-pseudolattice and $Idl(A)$ be
the frame of all ideals of $A$. If $J\in Idl(A)$ then we will
write  $\neg J$ for the pseudocomplement of
$J$ in $Idl(A)$ (i.e. $\neg J=\bv\{I\in Idl(A)\st I\we
J=\{0\}\}$). Note that $\neg J=\{a\in A\st (\fa b\in J)(a\we
b=0)\}$ (see Stone \cite{ST}). Recall that an ideal $J$ of $A$ is
called {\em simple} (Stone \cite{ST}) if $J\vee\neg J= A$ (i.e. $J$ has a complement in $Idl(A)$). As it
is proved in \cite{ST}, the set $Si(A)$ of all simple ideals of
$A$ is a Boolean algebra with respect to the lattice operations in
$Idl(A)$.
\end{nist}

\begin{fact}\label{gbpapi}{\rm (\cite{Di4})}
(a) A distributive $\{0\}$-pseudolattice $A$ is a generalized
Boolean pseudolattice iff every principal ideal of $A$ is simple.

\smallskip

\noindent(b) If $A$ is a generalized Boolean pseudolattice then
the correspondence $e_A:A\lra Si(A)$, $a\mapsto\dar(a)$, is a
dense $\{0\}$-pseudolattice embedding of $A$ in the Boolean
algebra $Si(A)$ and the pair $(Si(A),e_A(A))$ is an LBA.

\smallskip

\noindent(c){\rm (M. Stone \cite{ST})} An ideal of a Boolean
algebra is simple iff it is principal.
\end{fact}

\begin{fact}\label{gbplplba}{\rm (\cite{Di4})}
Let $(A,I)$ be an LBA. Then $I$ is a generalized Boolean pseudolattice.
\end{fact}

\begin{lm}\label{crhoslemma}{\rm (\cite{Di4})}
Let $(A,I)$ be an LBA and $\s\sbe A$. Then $\s$ is a bounded
cluster in $(A,\rho_s,I)$ iff it
 is a bounded ultrafilter in $A$.
\end{lm}

\begin{notas}\label{kxckx}
\rm Let $X$ be a topological space. We will denote by $CO(X)$ the
set of all clopen subsets of $X$,
and by $CK(X)$ the set of all clopen compact subsets of $X$. For
every $x\in X$, we set
$$ \s_x^{CO(X)}=\{F\in CO(X)\st x\in F\}
.$$
When there is no ambiguity, we will write $``\s_x^C$" instead of
$``\s_x^{CO(X)}$".
\end{notas}

\begin{defi}\label{defzlba}{\rm (\cite{Di4})}
\rm An LBA $(B, I)$  is called a {\em ZLB-algebra} (briefly, {\em
ZLBA}) if, for every $J\in Si(I)$, the join $\bv_B J$($=\bv_B
\{a\st a\in J\}$) exists.

Let $\ZLBA$ be the category whose objects are all ZLBAs and whose
morphisms are all functions $\p:(B, I)\lra(B_1, I_1)$ between the
objects of  $\ZLBA$ such that $\p:B\lra B_1$ is a Boolean
homomorphism satisfying the following condition:

\smallskip

\noindent(ZLBA) For every $b\in I_1$ there exists $a\in I$ such
that $b\le \p(a)$;

\smallskip

\noindent let the composition between the morphisms of $\ZLBA$ be
the usual composition between functions, and the
$\ZLBA$-identities be the identity functions.
\end{defi}

The Stone's result \ref{gbpapi}(c) 
leads to the following example:

\begin{exa}\label{zlbaexa}{\rm (\cite{Di4})}
\rm Let $B$ be a Boolean algebra. Then the pair $(B,B)$ is a ZLBA.
\end{exa}

\begin{theorem}\label{genstonec}{\rm (\cite{Di4})}
The category\/ $\ZLC$ of all locally compact zero-dimensio\-nal
Hausdorff spaces and all continuous maps between them is dually
equivalent to the category\/ $\ZLBA$. In more details, let $\TE^a:\ZLBA\lra\ZLC$ and
$\TE^t:\ZLC\lra\ZLBA$ be two contravariant functors defined as follows: for every
$X\in\card{\ZLC}$, $\TE^t(X)=(CO(X),  CK(X))$,  and for every $f\in\ZLC(X,Y)$,
$\TE^t(f):\TE^t(Y)\lra\TE^t(X)$ is defined by the formula
$\TE^t(f)(G)=f\inv(G), \ \  \fa G\in CO(Y)$;
 further, for every ZLBA $(B, I)$, $\TE^a(B, I)=\Psi^a(B,\rho_s,I)$
 (see \cite{Di3} for $\Psi^a$),  and
 for every $\p\in\ZLBA((B, I),(B_1, J))$,
 $\TE^a(\p):\TE^a(B_1, J)\lra\TE^a(B, I)$ is given by the formula
$\TE^a(\p)(u\ap)=\p\inv(u\ap), \ \  \fa u\ap\in\TE^a(B_1,J)$;
then $t^{C}:Id_{\ZLC}\lra   \TE^a\circ  \TE^t$,
where, $\fa X\in\card\ZLC$,  $t^{C}(X)=t_X^C$, and  $t_X^C(x)=\s_x^C$
for every $x\in X$, is a natural isomorphism; also,
$\l^C: Id_{\ZLBA}\lra  \TE^t\circ  \TE^a$, where, $\fa (B, I)\in\card\ZLBA$, $ \l^C(B,
I)=\l_B^C$, and $\l_B^C(a)=\l_B^g(a)$ for every $a\in B$ (see \cite{Di3} for $\l_B^g$), is a natural isomorphism.
\end{theorem}

\begin{cor}\label{stfact}{\rm (\cite{Di4})}
For every ZLBA $(B,I)$, the map $\SI_{(B,I)}:Si(I)\lra B$, where
$\SI_{(B,I)}(J)=\bv_B \{a\st a\in J\}$ for every $J\in Si(I)$, is
 a Boolean isomorphism.
\end{cor}

The next assertion was proved in \cite{Di4} although not stated explicitly there; it is a part of the proof of \cite[Theorem 1.14]{Di4}.

\begin{pro}\label{cox}{\rm (\cite{Di4})}
Let $(B,I)$ be a ZLBA and $X=\TE^a(B,I)(=\Psi^a(B,\rho_s,I))$. Then $\lbg(B)=CO(X)$ and $\lbg(I)=CK(X)$.
\end{pro}

\begin{defi}\label{defplbap}{\rm (\cite{Di4})}
\rm We will denote by $\PZLBA$  the cofull subcategory of the ca\-tegory $\ZLBA$
 whose
morphisms $\p:(A,I)\lra (B,J)$ satisfy the following additional
condition:

\smallskip

\noindent(PZLBA) $\p(I)\sbe J$.
\end{defi}

\begin{theorem}\label{genstonep}{\rm (\cite{Di4})}
The category\/ $\PZLC$ of all locally compact zero-dimensio\-nal
Hausdorff spaces and all perfect maps between them is dually
equivalent to the category\/ $\PZLBA$. In more details, the
restrictions $\TE^t_p:\PZLC\lra\PZLBA$ and
$\TE^a_p:\PZLBA\lra\PZLC$ of the duality functors $\TE^t$ and  $\TE^a$
defined in Theorem \ref{genstonec} are the
desired duality functors.
\end{theorem}

\begin{theorem}\label{genstonegstm}{\rm (\cite{Di4})}
The category\/ $\PZLBA$  is  equivalent to the category\/
$\GBPL$ whose objects are all generalized Boolean pseudolattices
and whose morphisms are all $\{0\}$-pseudolattice homomorphisms
between its objects satisfying condition (ZLBA). In more details, let
$E^a:\PZLBA\lra\GBPL$ and
$E^b:\GBPL\lra\PZLBA$ be two functors defined as follows: for every
$I\in\card{\GBPL}$, $E^b(I)=(Si(I),  e_I(I))$
(see \ref{sii} and \ref{gbpapi}(b) for the last notations),
and for every $\p\in\GBPL(I,J)$,
$E^b(\p):E^b(I)\lra E^b(J)$ is defined by the formula
$$E^b(\p)(I_1)=\bigcup\{\downarrow(\p(a))\st a\in I_1\}, \ \  \fa I_1\in Si(I);$$
 for every PZLBA $(B, I)$, $E^a(B, I)=I$,
   and
 for every $\p\in\PZLBA((B, I),(B_1, J))$,
 $E^a(\p):E^a(B, I)\lra E^a(B_1, J)$ is given by the formula
$E^a(\p)(a)=\p(a), \ \  \fa a\in E^a(B,I)$;
then $e^r:Id_{\GBPL}\lra   E^a\circ  E^b$,
where $e^r(I)=e_I^r, \ \ \fa I\in\card\GBPL$ and  $e_I^r:I\lra e_I(I)$
is the restriction of the map $e_I$ defined in \ref{gbpapi}(b), is a natural isomorphism; also,
$s: Id_{\PZLBA}\lra  E^b\circ  E^a$, where $ s(B,
I)=(\SI_{(B,I)})\inv, \ \ \fa (B, I)\in\card\PZLBA$ (see \ref{stfact} for $\SI_{(B,I)}$), is a natural isomorphism.

Hence, the categories\/ $\PZLC$ and $\GBPL$ are dually equivalent and the duality functors between them are $\TE^a_g:\GBPL\lra\PZLC$, $\TE^t_g:\PZLC\lra\GBPL$, where $\TE^a_g=\TE^a_p\circ E^b$ and $\TE^t_g=E^a\circ\TE^t_p$.
\end{theorem}

\section{Some New Duality Theorems}

Recall that
a homomorphism $\p$ between two Boolean algebras is called {\em complete}\/ if it
preserves all joins (and, consequently, all meets) that happen to exist;
this means that if $\{a_i\}$ is a family of elements in the domain of $\p$ with join
$a$, then the family $\{\p(a_i)\}$ has a join and that join is
equal to $\p(a)$.

\begin{defi}\label{skenuldim}
\rm
We will denote by $\SkeZLC$  the category of all zero-dimensional
locally compact Hausdorff spaces and all skeletal maps between them
(see \cite{Di2} for the definition of a skeletal map).

Let $\SkeZLBA$ be the cofull subcategory of the category $\ZLBA$  whose morphisms
are, in addition, complete homomorphisms.
\end{defi}

\begin{theorem}\label{skenuldimth}
The categories $\SkeZLC$ and $\SkeZLBA$ are dually equivalent.
\end{theorem}

\doc Having in mind Theorem \ref{genstonec}, it is enough to prove that if $f$ is a
$\SkeZLC$-morphism then $\TE^t(f)$ is  complete, and if $\p$ is a
$\SkeZLBA$-morphism
then $\TE^a(\p)$ is a skeletal map.

So, let $f\in\SkeZLC(X,Y)$ and $\p=\TE^t(f)$. Then
$\p:(CO(Y),CK(Y))\lra(CO(X),CK(X))$ and $\p(G)=f\inv(G)$, for all
$G\in CO(Y)$. Let $\{G_\g\st\g\in\GA\}\sbe CO(Y)$ and let this family have a join $G$ in $CO(Y)$. Set $W=\bigcup\{G_\g\st\g\in\GA\}$.
Since $Y$ is zero-dimensional, we get easily that $G=\cl(W)$. Thus
$\p(G)\spe\cl(\bigcup\{\p(G_\g)\st\g\in\GA\})=F$. Let $x\in f\inv(G)(=\p(G))$.
Then $f(x)\in G$ and there exists a
neighborhood $U$ of $x$ such that $f(U)\sbe G$. Suppose that $x\nin F$. Then there exists a neighborhood $V$ of $x$
such that $V\sbe U$ and $V\cap f\inv(G_\g)=\ems$ for all $\g\in\GA$. Thus $f(V)\cap W=\ems$. Then $\cl(f(V))\cap W=\ems$. Since
$\cl(f(V))\sbe\cl(f(U))\sbe G=\cl(W)$, we get that $\cl(f(V))\sbe \cl(W)\stm W(=\Fr(W))$. This leads to a contradiction because
$f$ is skeletal and thus $\int(\cl(f(V)))\nes$ (by \cite[Lemma 2.4]{Di2}).  So,
$\p(G)=f\inv(G)=F$. Since $\p(G)$ is clopen, we get that $\p(G)$ is the join of the family $\{\p(G_\g)\st\g\in\GA\}$ in $CO(X)$. Therefore, $\p$
is  complete.

Let now $\p\in\SkeZLBA((A,I),(B,J))$ and $f=\TE^a(\p)$.
Set $X=\TE^a(A,I)$ and $Y=\TE^a(B,J)$. Then
$f:Y\lra X$. Since $CK(Y)$ is an open base of $Y$,
for proving that $f$ is skeletal it is enough
to show that for every $G\in CK(Y)\stm\{\ems\}$, $\int(f(G))\nes$.
So, let $G\in CK(Y)\stm\{\ems\}$. Then there exists $b\in J\stm\{0\}$ such that
$G=\lbg(b)$. Suppose that $\bw\{c\in A\st b\le\p(c)\}=0$. Then, using the completeness of $\p$, we get
that $0=\p(0)=\bw\{\p(c)\st c\in A, b\le\p(c)\}\ge b$.
Since $b\neq 0$, we get a contradiction. Hence there exists $a\in A\stm\{0\}$ such that $a\le c$ for
all $c\in A$ for which $b\le\p(c)$.
 We will prove that $\lag(a)\sbe f(\lbg(b))(=f(G))$.
This will imply that $\int(f(G))\nes$. Let $u\in\lag(a)$. Then $a\in u$.
Suppose that there exists $c\in u$ such that $b\we\p(c)=0$. Then $b\le\p(c^*)$.
Thus $a\le c^*$, i.e. $a\we c=0$.
Since $a,c\in u$, we get a contradiction. Therefore, the set $\{b\}\cup\p(u)$ is a filter-base.
Hence there exists an ultrafilter $v$ in $B$ such that $\{b\}\cup\p(u)\sbe v$.
Then $b\in v$ and $u\sbe\p\inv(v)$. Thus $u=\p\inv(v)$, i.e. $f(v)=u$. So, $u\in f(\lbg(b))$. \sqs

\begin{rems}\label{equivcond}
\rm Note that in the definition of the category $\SkeZLBA$ the requirement that the
morphisms
$\p:(A,I)\lra(B,J)$ are complete can be replaced by the  following condition:

\smallskip

\noindent (SkeZLBA)  For every $ b\in J\stm\{0\}$ there exists $a\in I\stm\{0\}$ such that
$(\fa c\in A)[(b\le\p(c))\rightarrow(a\le c)]$.

\smallskip

\noindent Indeed, the proof of the above theorem shows the sufficiency of this condition and its necessity
can be established as follows. Let $f\in\SkeZLC(X,Y)$ and $\p=\TE^t(f)$. Then
$\p:(CO(Y),CK(Y))\lra(CO(X),CK(X))$ and $\p(G)=f\inv(G)$, for all
$G\in CO(Y)$. Let $F\in CK(X)\stm\{\ems\}$. Then
$\int (f(F))\nes$. Hence there exists $G\in CK(Y)\stm\{\ems\}$
such that $G\sbe\int(f(F))$. Let $H\in CO(Y)$ and $F\sbe f\inv(H)$. Then
$G\sbe\int(f(F))\sbe f(F)\sbe H$. So, condition (SkeZLBA) is satisfied.

Moreover, condition (SkeZLBA) can be replaced by the following one:

\smallskip

\noindent (CEP)  For every $ b\in B\stm\{0\}$ there exists $a\in A\stm\{0\}$ such that
$(\fa c\in A)[(b\le\p(c))\rightarrow(a\le c)]$.

\smallskip

\noindent Indeed, if $b\in B\stm\{0\}$ then, by the density of $J$ in $B$, there exists $b_1\in I\stm\{0\}$
such that $b_1\le b$. Now, applying (SkeZLBA) for $b_1$, we get an $a\in I\stm\{0\}$ which satisfies also the requirements of (CEP)
about $b$.
Conversely, if $ b\in J\stm\{0\}$ then, by (CEP), there exists  $a\in A\stm\{0\}$ such that
$(\fa c\in A)[(b\le\p(c))\rightarrow(a\le c)]$; but, by condition (ZLBA) (see \ref{defzlba}), there exists $a_1\in I$ such that $b\le\p(a_1)$;
thus $a\le a_1$; since $I$ is an ideal, we get that $a\in I$; so, condition (SkeZLBA) is satisfied.
\end{rems}

We denote by $S^t:\ZHC\lra\Bool$ and $S^a:\Bool\lra \ZHC$ the Stone duality functors
between the categories $\ZHC$ of compact zero-dimensional Hausdorff spaces
and continuous maps and  $\Bool$ of Boolean algebras and Boolean homomorphisms.

The assertion (c) of the next corollary is a zero-dimensional analogue of the Fedorchuk Duality Theorem \cite{F2}.

\begin{cor}\label{corperqo}
(a) Let $f$ be a $\PZLC$-morphism. Then $f$ is a quasi-open map (see \cite{Di2} for this notion)
iff $\TE^t(f)$ is complete.
In particular, if $f$ is a $\ZHC$-morphism then $f$ is a quasi-open map
iff $S^t(f)$ is complete.

\smallskip

\noindent(b) The cofull subcategory $\QPZLC$ of the category $\PZLC$ (see \ref{genstonep}) whose morphisms are, in addition, quasi-open maps, is dually equivalent to the cofull subcategory $\QPZLBA$
of the category $\PZLBA$ whose morphisms are, in addition, complete homomorphisms;

\smallskip

\noindent(c) The category $\QZHC$ of compact zero-dimensional Hausdorff spaces  and quasi-open maps
is dually equivalent to the category $\BoolC$
of  Boolean algebras and complete  Boolean homomorphisms.
\end{cor}

\doc The assertion (a) follows from the proof of Theorem \ref{skenuldimth} and \cite[Corollary 2.5]{Di2}.
The assertions (b) and (c) follow from (a) and  Theorem \ref{skenuldimth}. \sqs

The last corollary together with Fedorchuk Duality Theorem \cite{F2} imply the following assertion in which the equivalence $(a)\iff(b)$
is a special case of a
much more general theorem due to Monk \cite{Monk}.

\begin{cor}\label{cormonk}
Let $\p\in\Bool(A,B)$ and $A\ap,B\ap$ be  minimal completions of $A$ and $B$ respectively. We can suppose that $A\sbe A\ap$ and $B\sbe B\ap$. Then the following conditions are equivalent:

\smallskip

\noindent(a) $\p$ can be extended to a complete homomorphism $\psi:A\ap\lra B\ap$;

\smallskip

\noindent(b) $\p$ is a complete homomorphism;

\smallskip

\noindent(c) $\p$ satisfies condition (CEP) (see \ref{equivcond} above).
\end{cor}

\doc (a)$\Rightarrow$(b) This is obvious.

\noindent(b)$\Rightarrow$(c) This was already established in the proof of Theorem \ref{skenuldimth} (see also \ref{equivcond}).

\noindent(c)$\Rightarrow$(a) Obviously, $\p\in\ZLBA((A,A),(B,B))$. Then, by \ref{equivcond} and Theorem \ref{skenuldimth}, condition (CEP) implies that $f=\TE^a(\p)(=S^a(\p))$
is a skeletal map. Since $f$ is  closed, we get that $f$ is a quasi-open map between $Y=\TE^a(B,B)(=S^a(B))$ and $X=\TE^a(A,A)(=S^a(A))$.
Now, by Fedorchuk Duality Theorem \cite{F2}, the map $\psi=\Psi^t(f):RC(X)\lra RC(Y)$ is a complete homomorphism. Since, for every $F\in RC(X)$, $\psi(F)=\cl(f\inv(\int(F)))$ and $CO(X)$ is a subalgebra of the Boolean algebra $RC(X)$, we get that for every $F\in CO(X)$, $\psi(F)=f\inv(F)=\p\ap(F)$
(here $\p\ap=\TE^t(\TE^a(\p))$).  Then the existence of a natural isomorphism between  the composition $\TE^t\circ\TE^a$ and the identity functor (see Theorem \ref{genstonec}),
and the fact that $RC(X)$ and $RC(Y)$ are minimal completions of, respectively, $A$ and $B$, imply our assertion. \sqs

Now, using Theorem \ref{genstonegstm}, we will present in a simpler form the result established in Corollary \ref{corperqo}(b).

\begin{theorem}\label{thqgbpl}
The category $\QPZLC$ is dually equivalent to the cofull subcategory $\QGBPL$ of the category $\GBPL$ whose morphisms, in addition, preserve all meets that happen to exist.
\end{theorem}

\doc Having in mind Theorem \ref{genstonegstm} and Corollary \ref{corperqo}(b), it is enough to show that the functor $E^a$ (see \ref{genstonegstm}) maps $\QPZLBA$ to $\QGBPL$ and the functor $E^b$ (see again \ref{genstonegstm}) maps $\QGBPL$ to $\QPZLBA$ because with this we will obtain that the categories $\QPZLBA$ and $\QGBPL$ are equivalent. Obviously, if $\p\ap:(A,I)\lra (B,J)$ is a $\QPZLBA$-morphism then $\p=E^a(\p\ap)=(\p\ap)_{|I}:I\lra J$ preserves all meets in $I$ that happen to exist (indeed, since $I$ is an ideal of $A$, every meet in $I$ of elements of $I$ is also a meet in $A$). Conversely, let $\p:I\lra J$ be a $\QGBPL$-morphism. We will show that $\p$ satisfies the following condition:

\smallskip

\noindent(QGBPL) For every $ b\in J\stm\{0\}$ there exists $a\in I\stm\{0\}$ such that
$(\fa c\in I)[(b\le\p(c))\rightarrow(a\le c)]$.

\smallskip

\noindent Indeed, let $ b\in J\stm\{0\}$. Suppose that $\bw_I\{c\in I\st b\le\p(c)\}=0$. Then, using the completeness of $\p$, we get
that $0=\p(0)=\bw\{\p(c)\st c\in I, b\le\p(c)\}\ge b$.
Since $b\neq 0$, we get a contradiction. Hence there exists $a\in I\stm\{0\}$ such that $a\le c$ for
all $c\in I$ for which $b\le\p(c)$.

Let $\p\ap=E^b(\p)$. We will show that the map $\p\ap$ satisfies condition (SkeZLBA). We have that $\p\ap:(SI(I),e_I(I))\lra(Si(J),e_J(J))$. Let $J_1\in e_J(J)\stm\{0\}$.
Then there exists $b\in J\stm\{0\}$ such that $J_1=\downarrow(b)$. By (QGBPL), there exists $a\in I\stm\{0\}$ such that
$(\fa c\in I)[(b\le\p(c))\rightarrow(a\le c)]$. Let $I_1\in Si(I)$ and $J_1\sbe\p\ap(I_1)$. Then, by the definition of the map $\p\ap$ (see Theorem \ref{genstonegstm}),
we have that $\downarrow(b)\sbe\bigcup\{\downarrow(\p(c))\st c\in I_1\}$. Thus there exists $c\in I_1$ such that $b\le\p(c)$. Since $c\in I$, we get that $a\le c$.
Therefore, $\downarrow(a)\sbe I_1$. So, the map $\p\ap$ satisfies condition (SkeZLBA). Now \ref{equivcond} implies that $\p\ap$ is a complete homomorphism. Thus $\p\ap$ is a $\QPZLBA$-morphism. \sqs

\begin{rem}\label{equivcond1}
\rm The proof of Theorem \ref{thqgbpl} shows that in the definition of the category $\QPZLBA$ the requirement that its
morphisms
$\p:I\lra J$ preserve all meets that happen to exist can be replaced by the condition (QGBPL) introduced above.
\end{rem}

\begin{theorem}\label{thzop}
(a) Let $f\in\ZLC(X,Y)$, $\p=\TE^t(f)$, $(A,I)=\TE^t(X)$ and $(B,J)=\TE^t(Y)$. Then the map $f$ is open iff there exists a map $\psi:I\lra J$ which satisfies the following conditions:

\smallskip

\noindent{\rm (OZL1)} For every $b\in J$ and every $a\in I$, $(a\we\p(b)=0)\rightarrow (\psi(a)\we b=0)$, and

\smallskip

\noindent{\rm (OZL2)} For every $a\in I$, $\p(\psi(a))\ge a$

\smallskip

\noindent(such a map $\psi$ will be called a \/ {\em lower P-preadjoint} of $\p$).

\smallskip

\noindent(b) The cofull subcategory $\OZLC$ of the category $\ZLC$ whose morphisms are, in addition, open maps is dually equivalent to the cofull subcategory
$\OZLBA$ of the category $\ZLBA$ whose morphisms have, in addition, lower P-preadjoints.
\end{theorem}

\doc (a) Let $f\in\ZLC(X,Y)$ be an open map. For every $F\in CK(X)(=I)$ set $\psi(F)=f(F)$. Then, clearly, $\psi(F)\in CK(Y)(=J)$ and $\p(\psi(F))=f\inv(f(F))\spe F$.
Hence, condition (OZL2) is satisfied. Let $F\in CK(X)$, $G\in CK(Y)$  and $F\we\p(G)=0$. Then $F\cap f\inv(G)=\ems$. Thus $f(F)\cap G=\ems$, i.e. $\psi(F)\we G=0$.
Therefore, condition (OZL1) is satisfied as well.

Let now $\p$ has a lower P-preadjoint. We will show that $f\ap=\TE^a(\p)$ is an open map. This will imply that $f$ is open. Let $X\ap=\TE^a(\TE^t(X))$ and $Y\ap=\TE^a(\TE^t(Y))$. Since $\lag(I)$ is an open base of $X\ap$, it is enough to show that $f\ap(\lag(a))$ is an open set, for every $a\in I$. So, let $a\in I$. We will prove that $f\ap(\lag(a))=\lbg(\psi(a))$. Let $u\in\lag(a)$. Then $a\in u$. Let $v=f\ap(u)$, i.e. $v=\p\inv(u)$. By (OZL2), $\p(\psi(a))\ge a$ and hence   $\p(\psi(a))\in u$. Thus $\psi(a)\in\p\inv(u)=v$, i.e. $f\ap(u)\in\lbg(\psi(a))$. Therefore $f\ap(\lag(a))\sbe\lbg(\psi(a))$.
Conversely, let $v\in\lbg(\psi(a))$. Then $\psi(a)\in v$. Suppose that there exists $b\in v$ such that $a\we\p(b)=0$. Since $v$ is a bounded ultrafilter, there exists $b_0\in v\cap J$. Then $b_1=b\we b_0\in J\cap v$ and $a\we\p(b_1)=0$. Now, condition (OZL1) implies that $\psi(a)\we b_1=0$, which is a contradiction. Hence, the set $\{a\}\cup\p(v)$ is a filter-base. Thus there exists an ultrafilter $u\spe \{a\}\cup\p(v)$. Then $a\in u\cap I$ and $v\sbe\p\inv(u)$. Therefore, $v=\p\inv(u)=f(u)$. This shows that
$f\ap(\lag(a))\spe\lbg(\psi(a))$.  Hence, $f\ap$ is an open map.

\smallskip

\noindent(b) It follows from (a) and Theorem \ref{genstonec}. \sqs

\begin{rems}\label{remop2}
\rm
Note that condition (OZL2) implies condition (ZLBA). Indeed, in the notations of Theorem \ref{thzop}, if $a\in I$ then $b=\psi(a)\in J$ and $\p(b)\ge a$.
Further, condition (OZL1) implies that (again in the notations of Theorem \ref{thzop}) $\psi(0)=0$. Indeed,  $0\we\p(\psi(0))=0$ implies that $\psi(0)\we\psi(0)=0$, i.e. that $\psi(0)=0$.
\end{rems}

\begin{nist}\label{ladj}
\rm Let us recall the notion of {\em lower adjoint} for posets. Let $\p:A\lra B$ be an order-preserving map between posets.
If $\psi:B\lra A$ is an order-preserving map satisfying the following condition

\smallskip

 \noindent($\LAM$) for all $a\in A$ and all
$b\in B$, $b\le \p(a)$ iff $\psi (b)\le a$

\smallskip

\noindent(i.e. the pair
$(\psi,\p)$ forms a Galois connection between posets $B$ and $A$) then we will say that $\psi$ is a {\em lower adjoint}\/ of $\p$.
It is easy to see that condition ($\LAM$) is equivalent to  the following  condition:

\smallskip

\noindent($\LAM\ap$) $\fa a\in A$ and $\fa b\in B$, $\psi (\p(a))\le a$ and $\p(\psi (b))\ge b$.
\end{nist}

\begin{theorem}\label{thzopp}
(a) Let $f\in\PZLC(X,Y)$,  $(A,I)=\TE^t(X)$, $(B,J)=\TE^t(Y)$ and $\p=\TE^t(f)$. Then the map $f$ is open iff $\p:B\lra A$ has a lower adjoint $\psi:A\lra B$.

\smallskip

\noindent(b) The cofull subcategory $\OPZLC$ of the category $\PZLC$ whose morphisms are, in addition, open maps is dually equivalent to the cofull subcategory
$\OPZLBA$ of the category $\PZLBA$ whose morphisms have, in addition, lower adjoints.
\end{theorem}

\doc (a) Let $f\in\PZLC(X,Y)$ and $f$ is open. Then set $\psi(F)=f(F)$, for every $F\in CO(X)(=A)$. Then, since $f$ is open and closed map, $\psi:A\lra B(=CO(Y))$. Obviously, $\p(\psi(F))=f\inv(f(F))\spe F$ for every $F\in CO(X)$ and $\psi(\p(G))=f(f\inv(G))\sbe G$ for every $G\in CO(Y)$. Hence $\psi$ is a lower adjoint of $\p$. Conversely, let $\p$ has a lower adjoint $\psi$. Then $\psi(I)\sbe J$. Indeed, let $a\in I$. Then, by condition (ZLBA), there exists $b\in J$ such that $a\le \p(b)$. Then $\psi(a)\le\psi(\p(b))\le b\in J$. Thus, $\psi(a)\in J$. Further, condition (OZL2) is clearly fulfilled as well as condition (OZL1) (see, e.g., \cite[Fact 1.22(a)]{Di2}). So, $(\psi)_{|I}$ is a lower P-preadjoint of $\p$. Then, by Theorem \ref{thzop}(a), $f:X\lra Y$ is an open map.

\smallskip

\noindent(b) It follows from (a) and Theorem \ref{genstonep}. \sqs

\begin{cor}\label{thzopps}
(a) Let $f\in\ZHC(X,Y)$, $\p=S^t(f)$, $A=S^t(X)$ and $B=S^t(Y)$. Then the map $f$ is open iff $\p:B\lra A$ has a lower adjoint $\psi:A\lra B$.

\smallskip

\noindent(b) The category $\OZHC$ of compact zero-dimensional Hausdorff spaces  and open maps
is dually equivalent to the category $\OBool$
of  Boolean algebras and  Boolean homomorphisms having lower adjoints.
\end{cor}

\doc It follows immediately from Theorem \ref{thzopp}. \sqs

\begin{defi}\label{thzoppgd}
\rm
Let $\p\in\GBPL(J,I)$. If $\psi:I\lra J$ is a map which satisfies conditions (OZL1) and (OZL2) (see \ref{thzop})
then $\psi$ is called a {\em lower preadjoint} of $\p$.

Let $\OGBPL$ be the cofull subcategory
of the category $\GBPL$ whose morphisms have, in addition, lower preadjoints.
\end{defi}

\begin{cor}\label{thzoppg}
The categories $\OPZLC$ and\/ $\OGBPL$ are dually equivalent.
\end{cor}

\doc It follows from Theorems \ref{genstonegstm}, \ref{thzop} and  \ref{thzopp}. Indeed, it is enough to show that the categories $\OGBPL$ and $\OPZLBA$ are equivalent.
 From the proof of Theorem \ref{thzopp}, it follows that if $\p\ap$ is an  $\OPZLBA$-morphism then $\p=E^a(\p\ap)$ has a lower preadjoint. Conversely, if $\p$ is an $\OGBPL$-morphism then $\p\ap=E^b(\p)$ can be regarded as an extension of $\p$. This implies immediately that $\p\ap$ has a lower P-preadjoint. Now, Theorem \ref{thzop} implies that $f=\TE^a(\p\ap)$ is an open map. Thus, by Theorem \ref{thzopp}, $\p\ap$ has a lower adjoint.  \sqs

\section{Surjective and Injective Maps. Embeddings}

 In this section we will investigate the following problem: characterize the injective and surjective   morphisms of the categories $\ZLC$ and $\HLC$ and their subcategories, discussed in \cite{Di2, Di3, Di4} and in the previous section here,  by means of some  properties of their dual morphisms. Such a problem was regarded by M. Stone in \cite{ST}  for the morphisms of the category $\ZLC$ (for surjective morphisms) and of the category $\PZLC$ (for injective morphisms), and by de Vries in \cite{dV2} for the morphisms of the category $\HC$ of compact Hausdorff spaces and continuous maps. We will generalize their results. An analogous problem will be investigated for the homeomorphic embeddings, LCA-embeddings and the morphisms of the categories dual to the mentioned above. Our result about LCA-embeddings (see Theorem \ref{lcaemb}) generalizes a theorem of  Fedorchuk  \cite[Theorem 6]{F2}. Obviously, in all categories regarded here the injective morphisms are monomorphisms and surjective morphisms are epimorphisms. It is clear that surjectivity characterizes epimorphisms in the categories $\HC$ and $\ZHC$, and injectivity characterizes the monomorphisms in the categories $\HLC$, $\ZLC$, $\PLC$ and $\PZLC$. For the other categories regarded here I have no characterization of the monomorphisms and epimorphisms. The things are not simple as shows the following assertion:

 \begin{pro}\label{monohqo}
 The closed irreducible maps are monomorphisms in the category of all Hausdorff spaces and quasi-open maps between them.
 \end{pro}

 \doc Note first that every closed irreducible map is quasi-open (see \cite{P}). Let now $X,Y,Z$ be Hausdorff spaces, $f:Y\lra Z$ be a closed irreducible map  and $h,k:X\lra Y$ be two quasi-open maps such that $f\circ h=f\circ k$. Suppose that there exists $x\in X$ such that $h(x)\neq k(x)$. Then there exist disjoint open sets $U$ and $V$ in $Y$ such that $h(x)\in U$ and $k(x)\in V$. Let $W=h\inv(U)\cap k\inv(V)$. Then $x\in W$, $h(W)\sbe U$ and $k(W)\sbe V$. Set $O=\int(h(W))$. Since $h$ is quasi-open, we get that $O\nes$. Let $F=Y\stm O$. Then $F$ is a closed proper subset of $Y$. We will show that $f(F)=Y$ which will be a contradiction. We have that $f(O)=f(\int(h(W)))\sbe f(h(W))=f(k(W))\sbe f(V)\sbe f(F)$. Now, the surjectivity of $f$ implies that $f(F)=Y$. \sqs

We start with some simple observations.

\begin{pro}\label{injdl}
Let $f\in\HLC(X,Y)$, $(A,\rho,\BBBB)=\LAM^t(X)$, $(B,\eta,\BBBB\ap)=\LAM^t(Y)$ and $\p=\LAM^t(f)$. Then $\p$ is an injection $\iff$ $\p_{|\BBBB\ap}$ is an injection  $\iff$ $\cl_Y(f(X))=Y$.
\end{pro}

\doc We have that $\p: RC(Y)\lra RC(X)$.

Obviously, if $\p$ is an injection then $\p_{|\BBBB\ap}$ is an injection.

Let $\p_{|\BBBB\ap}$ be an injection, $G\in CR(Y)$ and $G\nes$. Then $\p(G)\nes$, i.e. $f\inv(\int(G))\nes$. This means that $f(X)\cap\int(G)\nes$. Thus $\cl(f(X))=Y$.

Finally, let $\cl(f(X))=Y$, $G,H\in RC(Y)$, $G\neq H$ and $\p(G)=\p(H)$. Then, by the continuity of $f$, $\cl_Y(f(\cl_X(f\inv(\int_Y(G)))))= \cl_Y(f((f\inv(\int_Y(G))))= \cl_Y(f(X)\cap\int_Y(G))=G$ and, analogously,  $\cl_Y(f(\cl_X(f\inv(\int_Y(H)))))=H$. Hence $G=H$, a contradiction. So, $\p$ is an injection. \sqs

A simplified version of the above proof shows that the following assertion takes place:

\begin{pro}\label{injdz}
Let $f\in\ZLC(X,Y)$, $(A,I)=\TE^t(X)$, $(B,J)=\TE^t(Y)$ and $\p=\TE^t(f)$. Then $\p$ is an injection $\iff$ $\p_{|J}$ is an injection $\iff$ $\cl_Y(f(X))=Y$.
\end{pro}

\begin{pro}\label{surdz}
Let $f\in\ZLC(X,Y)$, $\p=\TE^t(f)$, $(A,I)=\TE^t(X)$, $(B,J)=\TE^t(Y)$ and $\p(B)\spe I$ (or $\p(J)\spe I$). Then $f$ is an injection.
\end{pro}

\doc Suppose that there exist $x,y\in X$ such that $x\neq y$ and $f(x)=f(y)$. Then there exists $U\in CK(X)$ such that $x\in U\sbe X\stm\{y\}$. There exists $V\in CO(Y)$ with $\p(V)=U$, i.e. $f\inv(V)=U$. Then $f(U)=f(X)\cap V$ and hence $f\inv(f(U))=X\cap f\inv(V)=U$. Since $f(y)=f(x)\in f(U)$, we get that $y\in U$, a contradiction. Thus, $f$ is an injection. \sqs

We proceed  with the investigation of the problem in the category $\ZLC$ and its subcategories.

\begin{theorem}\label{inzdlc}
Let $f\in\ZLC(X,Y)$, $\p=\TE^t(f)$, $(A,I)=\TE^t(X)$ and $(B,J)=\TE^t(Y)$. Then $f$ is an injection iff $\p:(B,J)\lra(A,I)$ satisfies the following condition:

\smallskip

\noindent{\rm (InZLC)} For any $a,b\in I$ such that $a\we b=0$ there exists $a\ap,b\ap\in J$ with $a\ap\we b\ap=0$, $\p(a\ap)\ge a$ and $\p(b\ap)\ge b$.
\end{theorem}

\doc Let $f:X\lra Y$ be an injection. We have that $(A,I)=(CO(X),CK(X))$, $(B,J)=(CO(Y),CK(Y))$ and $\p:CO(Y)\lra CO(X)$, $G\mapsto f\inv(G)$. Let $F_1,F_2\in CK(X)$ and $F_1\cap F_2=\ems$. Since $f$ is a continuous injection, we get that $f(F_1)$ and $f(F_2)$ are disjoint compact subsets of $Y$. Using the fact that $CK(Y)$ is a base of $Y$, we get that there exist disjoint $G_1,G_2\in CK(Y)$ such that $f(F_i)\sbe G_i$, $i=1,2$. Then  $F_i\sbe f\inv(G_i)$, i.e. $F_i\sbe\p(G_i)$, $i=1,2$. Hence, $\p$ satisfies condition (InZLC).

Let now  $\p$ satisfies condition (InZLC). We will prove that $f$ is an injection. Let $x,y\in X$ and $x\neq y$. Then there exist disjoint $F_x,F_y\in CR(X)$ such that $x\in F_x$ and $y\in F_y$. Now, by condition (InZLC), there exist $G_x,G_y\in CR(Y)$ such that $G_x\cap G_y=\ems$, $f\inv(G_x)\spe F_x$ and $f\inv(G_y)\spe F_y$. Then $f(x)\in G_x$ and $f(y)\in G_y$. Thus $f(x)\neq f(y)$. \sqs

\begin{cor}\label{inzdlccor}
The cofull subcategory $\InZLC$ of the category $\ZLC$ whose morphisms are, in addition, injective maps, is dually equivalent to the cofull subcategory $\DInZLC$ of the category $\ZLBA$ whose morphism satisfy, in addition, condition (InZLC).
\end{cor}

\doc It follows from Theorems \ref{inzdlc} and \ref{genstonec}. \sqs

{\em In the sequel, we will not formulate corollaries like that because they follow directly from the respective characterization of injectivity or surjectivity and the corresponding duality theorem.}


\begin{rem}\label{rem3534}
\rm
Let us show how Theorem \ref{inzdlc} implies Proposition \ref{injdz}. So, let $\p(B)\spe I$. Then $\p(J)\spe I$. Indeed, let $a\in I$; then, by condition (ZLBA), there exists $b_1\in J$ such that $\p(b_1)\ge a$; since there exists $b_2\in B$ with $\p(b_2)=a$, we get that $\p(b_1\we b_2)=a$ and $b_1\we b_2\in J$. Hence, $\p(J)\spe I$. Let now $a,b\in I$ and $a\we b=0$. There exist $a_1,b_1\in J$ such that $\p(a_1)=a$ and $\p(b_1)=b$. Then $\p(a_1\we b_1^*)=a\we b^*=a$, $a_1\we b_1^*\in J$ and $(a_1\we b_1^*)\we b_1=0$. Therefore, $\p$ satisfies condition (InZLC).
\end{rem}

In the next theorem we will assume that the ideals and prime ideals could be non-proper.

\begin{theorem}\label{surzdlc}
Let $f\in\ZLC(X,Y)$, $\p=\TE^t(f)$, $(A,I)=\TE^t(X)$ and $(B,J)=\TE^t(Y)$. Then the following conditions are equivalent:

\smallskip

\noindent(a) $f$ is a surjection;

\smallskip

\noindent(b) $\p:(B,J)\lra(A,I)$ is an injection and for every bounded ultrafilter $v$ in $(B,J)$ there exists $a\in I$ such that $a\we\p(v)\neq 0$ (i.e.
$a\we\p(b)\neq 0$ for any $b\in v$);

\smallskip

\noindent(c) $\p:(B,J)\lra(A,I)$ is an injection and for every prime ideal $J_1$ of $J$, we have that $\bv\{I_{\p(b)}\st b\in J_1\}=I$ implies $J_1=J$ (where $I_{\p(b)}=\{a\in I\st a\le\p(b)\}$);

\smallskip

\noindent(d) $\p:(B,J)\lra(A,I)$ is an injection and for every  ideal $J_1$ of $J$, $[(\bv\{I_{\p(b)}\st b\in J_1\}=I)\rightarrow(J_1=J)]$.
\end{theorem}

\doc (a)$\Rightarrow$(b) Let $f(X)=Y$. Then, for every $G\in CO(Y)$, $f\inv(G)=\ems$ iff $G=\ems$. This means that for every $G\in CO(Y)$, $\p(G)=0$ iff $G=0$. Hence, $\p$ is an injection. Further, the bounded ultrafilters in $(B,J)=(CO(Y),CK(Y))$ are of the form $\s_y^C$ (see \ref{kxckx} for this notation) and analogously for $(A,I)$. So, let $y\in Y$. Then there exists $x\in X$ such that $f(x)=y$. This implies that $\p(\s_y^C)\sbe\s_x^C$. There exists $F\in CK(X)\cap\s_x^C$. Then $F\cap f\inv(G)\nes$, for every $G\in\s_y^C$, i.e. $F\we\p(\s_y^C)\neq 0$.

\smallskip

\noindent(b)$\Rightarrow$(c) Let $J_1$ be a prime ideal of $J$.  Let $\bv\{I_{\p(b)}\st b\in J_1\}=I$. Suppose that $J_1\neq J$. Then $v_1=\{b\in B\st b\we(J\stm J_1)\neq 0\}$ is a bounded ultrafilter in $(B,J)$ and $v_1\cap J=J\stm J_1$ (indeed, apply \cite[Proposition 3.6]{Di3} to the LCA $(B,\rho_s,J)$). By (b), there exists $a\in I$ such that $a\we\p(v_1)\neq 0$. Since $a\in I$, there exist $b_1,\dots,b_k\in J_1$ and $a_1,\dots,a_k\in I$ (where $k\in\mathbb{N}^+$) such that $a=\bv\{a_i\st i=1,\dots,k\}$ and $a_i\le\p(b_i)$, $i=1,\dots,k$.  Set $b=\bv\{b_i\st i=1,\dots,k\}$. Then $a\le\p(b)$ and $b\in J_1$. Since $\p$ is an injection, we have that $\p(v_1\cap J)=\p(J\stm J_1)=\p(J)\stm\p(J_1)$. Thus $\p(b)\nin\p(v_1\cap J)$ (because $b\in J_1$). Since $a\le \p(b)$, we get that $\p(b)\we\p(v_1)\neq 0$. The injectivity of $\p$ implies that $b\we v_1\neq 0$. Thus $b\in v_1\cap J_1$. Therefore, $\p(b)\in \p(v_1\cap J_1)$, a contradiction. Hence, $J_1=J$.

\smallskip

\noindent(c)$\Rightarrow$(a) Suppose that $f(X)\neq Y$. Then there exists $y\in Y\stm f(X)$. Set $U=Y\stm\{y\}$. Thus $f(X)\sbe U$. Set $J_1=\{G\in CK(Y)\st G\sbe U\}$. Then $J_1$ is a prime ideal of $J(=CK(Y))$. (Indeed, if $G_1,G_2\in CK(Y)$ and $y\nin G_1\cap G_2$ then either $y\nin G_1$ or $y\in G_2$; hence, $G_1\in J_1$ or $G_2\in J_1$.) Obviously, $J_1\neq J$. We will prove that $\bv\{I_{\p(b)}\st b\in J_1\}=I$, which, by (c), will lead to a contradiction. So, let $F\in CK(X)$. Then $f(F)\sbe U$. Since $f(F)$ is compact, there exists $G\in CK(Y)$ such that $f(F)\sbe G\sbe U$. Then $G\in J_1$ and $F\sbe f\inv(G)=\p(G)$. Thus $F\in I_{\p(G)}$. Therefore,  $\bv\{I_{\p(b)}\st b\in J_1\}=I$. So, $f(X)=Y$.

\smallskip

\noindent(a)$\Rightarrow$(d) Let $f(X)=Y$. Then, as we have already proved (see (a)$\Rightarrow$(b)), $\p$ is an injection. Let $J_1$ be an ideal of $J$ such that $\bv\{I_{\p(b)}\st b\in J_1\}=I$. Suppose that $J_1\neq J$. Set $U=\bigcup\{G\st G\in J_1\}$. Then $U\neq Y$. (Indeed, if $U=Y$ then  every $H\in CK(Y)(=J)$ will be covered by a finite number of elements of $J_1$; since $J_1$ is an ideal, we will get that $H\in J_1$.) Since $f$ is a surjection, we get that $V=f\inv(U)\neq X$. Set $I_V=\{F\in I\st F\sbe V\}$. Then, obviously, $I_V$ is a proper ideal of $I$. Let $G\in J_1$ and $F\in I_{\p(G)}$. Then $F\sbe\p(G)= f\inv(G)\sbe f\inv(U)=V$. Thus
$\bv\{I_{\p(b)}\st b\in J_1\}\sbe I_V$. Since $I_V\neq I$, we get a contradiction. Therefore, $J_1=J$.

\smallskip

\noindent(d)$\Rightarrow$(c) It is obvious. \sqs

\begin{rem}\label{Stonesurj}
\rm
In \cite[Theorem 7]{ST} M. Stone proved a result which is equivalent to our assertion that
(a)$\Leftrightarrow$(d) in the previous theorem.  More precisely, M. Stone proved the following (in our notations): the map $f$ is a surjection iff the map $\psi=\p_{|J}:J\lra A$ is a 0-pseudolattice monomorphism and for every  ideal $J_1$ of $J$, $[(\bv\{I_{\p(b)}\st b\in J_1\}=I)\leftrightarrow(J_1=J)]$.
The Stone's condition $``(J_1=J)\rightarrow(\bv\{I_{\p(b)}\st b\in J_1\}=I)$", i.e. $``\bv\{I_{\p(b)}\st b\in J\}=I$", is equivalent (as it is easy to see) to our condition (ZLBA) (see \ref{defzlba}) which is automatically satisfied by the morphisms of the category $\ZLBA$ and thus it appears  in our Theorem \ref{surzdlc} in another form. Further, when $\p$ is an injection then, obviously, $\psi=\p_{|J}$ is an injection; in the converse direction we have the following: the map $\psi$ can be extended to a homomorphism $\p:B\lra A$ (by the result proved below) and then $\p$ is obliged to be an injection (indeed, if $b\in B\stm\{0\}$ and $\p(b)=0$ then the density of $J$ in $B$ implies that there exists $c\in J\stm\{0\}$ such that $c\le b$; then $\psi(c)=\p(c)=0$, a contradiction). So, our condition (d) is equivalent to the cited above Stone condition from \cite[Theorem 7]{ST}.
\end{rem}

\begin{pro}\label{exthom}
Let $(A,I)$ and $(B,J)$ be ZLBAs and $\psi:J\lra A$ be a 0-pseudolat\-tice homomorphism satisfying condition (ZLBA). Then $\psi$ can be extended to a homomorphic map $\p:B\lra A$.
\end{pro}

\doc For every $a\in A$ and every $b\in B$  set $I_a=\{c\in I\st c\le a\}$ and $J_b=\{c\in J\st c\le b\}$. It is easy to see that $I_a$ and $J_b$ are simple ideals of $I$ and $J$ respectively. Note also that $\neg I_a=I_{a^*}$ and analogously for $J_b$.

Let $b\in B$. Since $J$ is dense in $B$, we have that $b=\bv J_b$. We will show that $I(b)=\bv\{I_{\psi(c)}\st c\in J_b\}$ is a simple ideal of $I$. It is easy to see that $I(b)=\bigcup\{I_{\psi(c)}\st c\in J_b\}$. Let now $a\in I$. Then, by condition (ZLBA), there exists $c\in J$ such that $a\le\psi(c)$. We have that $c=(c\we b)\vee(c\we b^*)$, $c_1=c\we b\in J_b$, $c_2=c\we b^*\in\neg J_b$ and $c=c_1\vee c_2$. Thus $a\le\psi(c)=\psi(c_1)\vee\psi(c_2)$. We obtain that $a=a_1\vee a_2$, where $a_1=a\we\psi(c_1)$ and $a_2=a\we\psi(c_2)$. Obviously, $a_1\in I(b)$. We will show that $a_2\in\neg I(b)$. Indeed, let $a\ap\in I(b)$; then there exists $d\in J_b$ such that $a\ap\le\psi(d)$. Since $c_2\in\neg J_b$, we get that $d\we c_2=0$. Thus $\psi(d)\we\psi(c_2)=0$. Hence $a\ap\we a_2\le\psi(d)\we a\we\psi(c_2)=0$. Therefore, for every $a\ap\in I(b)$ we have that $a_2\we a\ap=0$. This means that $a_2\in\neg I(b)$. Therefore, $I(b)\vee \neg I(b)=I$, i.e. $I(b)$ is a simple ideal. Since $(A,I)$ is a ZLBA, we get that $\bv I(b)$ exists in $A$. We set now $\p(b)=\bv I(b)$. Obviously, $\p(0)=0$. Further, $\p(1)=\bv I(1)$. We have that $I(1)=\bigcup\{I_{\psi(c)}\st c\in J\}$. Applying condition (ZLBA), we get that $I(1)=I$. Now, using the density of $I$ in $A$, we obtain that $\p(1)=1$. Finally, the fact that $\p$ preserves finite meets and joins can be easily proved. Hence $\p:B\lra A$ is a Boolean homomorphism and the definition of $\p$ together with the density of $I$ in $A$ imply that $\p$ extends $\psi$. \sqs.

\begin{rem}\label{Stonesurjcor}
\rm
Note that \ref{Stonesurj} and \ref{exthom} imply that in Theorem \ref{surzdlc} we can obtain new conditions equivalent to the condition (a) by replacing in (b), (c) and (d) the phrase $``\p$ is an injection" by the phrase  $``\p_{|J}$ is an injection".
\end{rem}

\begin{theorem}\label{injopz}
Let $f\in\OZLC(X,Y)$, $\p=\TE^t(f)$, $(A,I)=\TE^t(X)$ and $(B,J)=\TE^t(Y)$. Then  $f$ is an injection $\iff$ $\p(J)\spe I$ $\iff$ $\p(B)\spe I$.
\end{theorem}

\doc  Note that, by Remark \ref{rem3534}, conditions $``\p(J)\spe I$" and $``\p(B)\spe I$" are equivalent.

 Let $f$ be an injection and $F\in CK(X)$. Then $f(F)\in CK(Y)$ and $f\inv(f(F))=F$. Hence, $\p(J)\spe I$.
Conversely, let $\p(J)\spe I$. Then, by \ref{surdz}, we get that $f$ is an injection. \sqs

\begin{theorem}\label{surperz}
Let $f\in\PZLC(X,Y)$, $\p=\TE^t(f)$, $(A,I)=\TE^t(X)$ and $(B,J)=\TE^t(Y)$. Then  $f$ is a surjection $\iff$ $\p$ is an injection $\iff$ $\p_{|J}$ is an injection.
\end{theorem}

\doc Clearly, Proposition \ref{injdz} implies that if $f$ is a surjection then  $\p$ is an injection. Hence $\p_{|J}$ is an injection.

Let now $\p_{|J}$ be an injection. Then, by Proposition \ref{injdz}, $\cl(f(X))=Y$. Since $f$ is a closed map, we get that $f$ is a surjection. \sqs

\begin{theorem}\label{injperz}
Let $f\in\PZLC(X,Y)$, $\p=\TE^t(f)$, $(A,I)=\TE^t(X)$ and $(B,J)=\TE^t(Y)$. Then  $f$ is an injection  iff $\p(J)=I$.
\end{theorem}

\doc Let $f$ be an injection. Then $f_{\upharpoonright X}:X\lra f(X)$ is a homeomorphism. Let $F\ap\in CK(X)$. Then $F=f(F\ap)$ is compact. Since $F$ is open in $f(X)$, there exists an open set $U$ in $Y$ such that $U\cap f(X)=F$. Then there exists $G\in CK(Y)$ such that $F\sbe G\sbe U$. Then, clearly, $f\inv(G)=f\inv(F)=F\ap $. Hence $\p(G)=F\ap$. Therefore, $\p(J)\spe I$. Since $f$ is perfect, we have that $\p(J)\sbe I$. Thus $\p(J)=I$. Conversely, let  $\p(J)=I$. Then Proposition \ref{surdz} implies that $f$ is an injection. \sqs

Obviously, the last two theorems imply the well-known Stone's results that a $\ZHC$-morphism $f$ is an injection (resp., a surjection) iff $\p=S^t(f)$ is a surjection (resp., an injection).

Now we turn to the morphisms of the category $\HLC$ and its subcategories  discussed in \cite{Di2,Di3}.


\begin{theorem}\label{surhlc}
Let $f\in\HLC(X,Y)$,  $(A,\rho,\BBBB)=\LAM^t(X)$, $(B,\eta,\BBBB\ap)=\LAM^t(Y)$ and $\p=\LAM^t(f)$. Then $f$ is a surjection iff $\p:(B,\eta,\BBBB\ap)\lra(A,\rho,\BBBB)$ satisfies the following condition:

\smallskip

\noindent{\rm (SuHLC)} For any bounded ultrafilter $v$ in $(B,\eta,\BBBB\ap)$ there exists a bounded ultrafilter $u$ in $(A,\rho,\BBBB)$ such that $\fa b\in\BBBB\ap$, $(b\eta v)\leftrightarrow((\fa b\ap\in \BBBB\ap)[(b\lle b\ap)\rightarrow(\p(b\ap)\rho u)])$.
\end{theorem}

\doc Let $f$ be a surjection and $v$ be a bounded ultrafilter  in $(B,\eta,\BBBB\ap)$. Recall that that  $(B,\eta,\BBBB\ap)=(RC(Y),\rho_Y,CR(Y))$. Obviously, $\bigcap v$ is an one-point set. Let $\{y\}=\bigcap v$. Since $f$ is a surjection, there exists $x\in X$ such that $f(x)=y$. Let $\s_x=\{F\in RC(X)\st x\in F\}$ and $\nu_x=\{F\in RC(X)\st
x\in\int(F)\}$. There exists an ultrafilter $u$ in $(A,\rho,\BBBB)$ such that $u\spe\nu_x$. Then it is easy to see that $u\sbe\s_x$. Hence $\bigcap u=\{x\}$. Let now $G,H\in CR(Y)$. Clearly, if $y\in G$ then $G\rho_Y v$ (i.e. $G\eta v$). Conversely, let $G\rho_Y v$. If $y\nin G$ then, by \cite[Corollary 3.1.5]{E2}, there exists $G\ap\in v$ such that $G\cap G\ap=\ems$, a contradiction. Hence, $G\rho_Y v$ iff $y\in G$. In an analogous way we obtain that $\cl(f\inv(\int(H)))\rho_X u$ iff $x\in\cl(f\inv(\int(H)))$. So, we have to show that $y\in G$ iff for all $H\in CR(Y)$, $(G\sbe\int(H))\rightarrow (x\in\cl(f\inv(\int(H))))$. Let $y\in G$ and $G\sbe\int(H)$ for some $H\in CR(Y)$. Then $\cl(f\inv(\int(H)))\spe f\inv(G)\spe f\inv(y)$ and thus $x\in\cl(f\inv(\int(H)))$. Conversely, let $x\in\cl(f\inv(\int(H)))$ for every $H\in CR(Y)$ such that $G\sbe\int(H)$. Suppose that $y\nin G$. Then there exists $H\in CR(Y)$ such that $G\sbe\int(H)\sbe H\sbe Y\stm\{y\}$. Then  $\cl(f\inv(\int(H)))\sbe f\inv(H)\sbe X\stm f\inv(y)$. Thus $x\nin\cl(f\inv(\int(H)))$, a contradiction. So, $\p$ satisfies condition (SuHLC).

Let now $\p$ satisfies condition (SuHLC). We will show that $f\ap=\LAM^a(\p)$ is a surjection. This will imply that $f$ is a surjection. Let $\s\ap\in Y\ap=\LAM^a(B,\eta,\BBBB\ap)$. Then there exists a bounded ultrafilter $v$ in $(B,\eta,\BBBB\ap)$ such that $\s\ap=\s_v$ (see \cite{Di2}). By condition (SuHLC),
there exists a bounded ultrafilter $u$ in $(A,\rho,\BBBB)$ such that $\fa b\in\BBBB\ap$, $(b\eta v)\leftrightarrow((\fa b\ap\in \BBBB\ap)[(b\lle b\ap)\rightarrow(\p(b\ap)\rho u)])$. Then it is easy to see that $f\ap(\s_u)=\s_v$ (see Theorem \ref{lccont} for the definition of the map $f\ap$). Hence, $f\ap$ is a surjection. \sqs

\begin{theorem}\label{injhlc}
Let $f\in\HLC(X,Y)$,  $(A,\rho,\BBBB)=\LAM^t(X)$, $(B,\eta,\BBBB\ap)=\LAM^t(Y)$ and $\p=\LAM^t(f)$. Then $f$ is an injection iff $\p:(B,\eta,\BBBB\ap)\lra(A,\rho,\BBBB)$ satisfies the following condition:

\smallskip

\noindent{\rm (InHLC)} For any $a,b\in\BBBB$, $a(-\rho) b$ implies that there exist $c,d\in\BBBB\ap$ such that $c\lle d$, $a\le\p(c)$ and $\p(d)(-\rho) b$.
\end{theorem}

\doc Let $f$ be an injection. We will show that $\p$ satisfies condition (InHLC). Let $F,G\in CR(X)$ and $F\cap G=\ems$. Since $f$ is an injection, we get that $f(F)\cap f(G)=\ems$. Using the fact that $f(F)$ and $f(G)$ are compact sets, we get that there exist $F\ap,G\ap\in CR(Y)$ such that $f(F)\sbe\int(F\ap)\sbe F\ap\sbe\int(G\ap)\sbe G\ap\sbe Y\stm f(G)$. Then, clearly, $F\sbe f\inv(\int(F\ap))\sbe\p(F\ap)$ and $G\cap\p(G\ap)=\ems$ (because $\p(G\ap)\sbe f\inv(G\ap)$ and $f\inv(G\ap)\cap G=\ems$). Therefore, $\p$ satisfies condition (InHLC).

Let now $\p$ satisfies condition (InHLC). We will prove that $f$ is an injection. Let $x,y\in X$ and $x\neq y$. Then there exist disjoint $F_x,F_y\in CR(X)$ such that $x\in F_x$ and $y\in F_y$. Now, by condition (InHLC), there exist $G_x,G_y\in CR(Y)$ such that $G_x\sbe\int(G_y)$, $F_x\sbe\cl(f\inv(\int(G_x)))$ and $\cl(f\inv(\int(G_y)))\cap F_y=\ems$.  Since $F_x\sbe f\inv(G_x)$, we get that $f(x)\in G_x$. Further, we have that $F_y\cap f\inv(\int(G_y))=\ems$. Thus $f(F_y)\cap\int(G_y)=\ems$. Then $f(F_y)\cap G_x=\ems$, and therefore,  $f(x)\neq f(y)$. Hence, $f$ is an injection. \sqs

The next theorem was proved in \cite[Theorem 2.11]{Di5} using a theorem of de Vries \cite{dV2}.  We will give now a direct proof of it
using only  \ref{injdl}.

\begin{theorem}\label{surplc}
Let $f\in\PLC(X,Y)$ and $\p=\LAM^t(f)$.  Then $f$ is a surjection iff $\p$ is an injection.
\end{theorem}

\doc If $f$ is a surjection then Proposition \ref{injdl} implies that $\p$ is an injection.
Let now $\p$ be an injection. Then, by Proposition \ref{injdl}, $\cl(f(X))=Y$. Since $f$ is a closed map, we get that $f$ is a surjection. \sqs

Obviously, Theorem \ref{injhlc} implies the following result:

\begin{theorem}\label{injplc}
Let $f\in\PLC(X,Y)$, $(A,\rho,\BBBB)=\LAM^t(X)$, $(B,\eta,\BBBB\ap)=\LAM^t(Y)$ and $\p=\LAM^t(f)$. Then $f$ is an injection iff $\p$ satisfies  condition
 (InHLC) (see Theorem \ref{injhlc}).
\end{theorem}

Another  characterization of the injective $\PLC$-morphisms was given  in \cite[Theorem 2.15]{Di5}. It was derived from a theorem of de Vries \cite{dV2}. It seems that the new condition looks better. In this form the theorem is new even in the compact case.

\begin{theorem}\label{injspzlc}
Let $f\in\PZLC(X,Y)$, $f$ be a skeletal map, $\p=\LAM^t(f)$ and $f$ be an injection. Then
$\p$ is a surjection.
\end{theorem}

\doc Let $\p\ap=\TE^t(f)$. Then, by Theorem \ref{injperz}, $\p\ap(CK(Y))= CK(X)$. Since $f$ is a skeletal map, \cite[Theorem 2.11]{Di2} (see also \ref{mainfed}) implies that $\p$ is a complete homomorphism.
Let $F\in RC(X)$. Then, using the fact that $CK(X)$ is a base of $X$, we get that $\int(F)=\bigcup\{G\in CK(X)\st G\sbe\int(F)\}$. Thus $F=\bv_{RC(X)}\{G\in CK(X)\st G\sbe\int(F)\}$. Further, for every
 $G\in CK(X)$ there exists $H_G\in CK(Y)$ such that $G=\p\ap(H_G)=\p(H_G)$. Therefore $F=\p(\bv_{RC(Y)}\{H_G\st G\in CK(X), G\sbe\int(F)\})$. This shows that $\p$ is a surjection. \sqs

\begin{rem}\label{monkinjsurj}
\rm
It is well known that the following two assertions can be add in Corollary \ref{cormonk}: $\p$ is an injection iff $\psi$ is an injection, and also, if $\p$ is a surjection then $\psi$ is a surjection (see, e.g., \cite[Ch. 25, Ex. 12]{GH}). Note that the first assertion follows immediately from Theorems \ref{surperz}   and \ref{surplc};
the second follows from  Theorems \ref{injperz} and \ref{injspzlc}.
\end{rem}

\begin{theorem}\label{slcompnk}
Let $f\in\SKLC(X,Y)$,  $(A,\rho,\BBBB)=\LAM^t(X)(=\Psi^t(X))$, $(B,\eta,\BBBB\ap)=\LAM^t(Y)(=\Psi^t(Y))$ and $\p=\LAM^t(f)(=\Psi^t(f))$. Then $f$ is a surjection if and only if $\p:(B,\eta,\BBBB\ap)\lra(A,\rho,\BBBB)$ satisfies the following condition:

\smallskip

\noindent{\rm (SuSkeLC)} For every bounded ultrafilter $v$ in
$(B,\eta,\BBBB\ap)$ there exists a bounded ultrafilter $u$ in
$(A,\rho,\BBBB)$ such that $\p\inv(u)\eta v$.
\end{theorem}

\doc Let $f:X\lra Y$ be a surjective continuous skeletal map
between two locally compact Hausdorff spaces and $\p=\Psi^t(f)$.
Then $\p:RC(Y)\lra RC(X)$ and $\pl(F)=\cl(f(F))$, for every $F\in
RC(X)$ (see the proof of \cite[Theorem 2.11]{Di2}). Let $v$ be a
bounded ultrafilter in $RC(Y)$. Then there exists $G_0\in
CR(Y)\cap v$. Hence there exists $y\in\bigcap\{G\st G\in v\}$.
Since $f$ is a surjection, there exists $x\in X$ such that
$f(x)=y$. Let $u$ be an ultrafilter in $RC(X)$ which contains
$\nu_x$ (see \cite[(3)]{Di2} for $\nu_x$). Then, obviously, $u$ is
a bounded ultrafilter in $(RC(X),\rho_X,CR(X))$. By
\cite[(51)]{Di2}, $u\sbe\s_x$ (see \cite[(3)]{Di2} for $\s_x$).
Hence $y\in\pl(F)$, for every $F\in u$. This means that
$\pl(u)\rho_Y v$. Since $\pl(u)$ is a filter-base of $\p\inv(u)$ (see \cite[(36)]{Di2}),
we get that $\p\inv(u)\rho_Y v$.
Therefore, $\p$ satisfies condition (SuSkeLC).

Let $\p$ satisfies condition (SuSkeLC). Set $f\ap=\Psi^a(\p)$. Let
$X\ap=\Psi^a(A,\rho,\BBBB)$, $Y\ap=\Psi^a (B,\eta,\BBBB\ap)$ and $\s\in
Y\ap$. Then $\s$ is a bounded cluster in $(B,\eta,\BBBB\ap)$. Hence there
exists a bounded ultrafilter $v$ in $(B,\eta,\BBBB\ap)$ such that
$\s=\s_v$. By (SuSkeLC), there exists a bounded ultrafilter $u$ in
$(A,\rho,\BBBB)$ such that $\p\inv(u)\eta v$. Thus $\p\inv(u)C_\eta
v$. Therefore, by  \cite[(35)]{Di2}, $f\ap(\s_u)=\s_{\pl(u)}=\s_v=\s$. So, $f\ap$ is a
surjection. Then, by  \cite[Theorem 2.11]{Di2}, $f$ is also a surjection.
 \sqs

Obviously, in condition (SuSkeLC), $``\p\inv(u)\eta v$" can be replaced by $``\pl(u)\eta v$".

Note that it is easy to see that condition (SuSkeLC) implies condition (SuHLC) when $\p$ is a $\SKAL$-morphism. This provides with a new proof  the sufficiency part of Theorem \ref{slcompnk}.

\begin{theorem}\label{ilcompnk}
Let $f\in\SKLC(X,Y)$,  $(A,\rho,\BBBB)=\LAM^t(X)(=\Psi^t(X))$, $(B,\eta,\BBBB\ap)=\LAM^t(Y)(=\Psi^t(Y))$ and $\p=\LAM^t(f)(=\Psi^t(f))$. Then $f$ is an injection if and only if $\p:(B,\eta,\BBBB\ap)\lra(A,\rho,\BBBB)$ satisfies the following condition:

\smallskip

\noindent{\rm (InSkeLC)} $\fa a,b\in \BBBB$, $\pl(a)\eta\pl(b)$ implies
$a\rho b$ (here $\pl$ is the left adjoint of $\p$).
\end{theorem}

\doc Let $f:X\lra Y$ be an injective continuous skeletal map.
The
function $\pl:RC(X)\lra RC(Y)$ is defined by $\pl(F)=\cl(f(F))$,
for every $F\in RC(X)$ (see \cite[(32) and (33)]{Di2}). Hence, for
$F\in CR(X)$, $\pl(F)=f(F)$. Since $f$ is an injection, it becomes
obvious that $\p$ satisfies condition (InSkeLC).

Let $\p$ satisfies condition (InSkeLC). We will show that $f$ is an injection. Let $x,y\in X$ and $x\neq y$. Then there exist disjoint $F_x,F_y\in CR(X)$ such that $x\in F_x$ and $y\in F_y$. If $f(x)=f(y)$ then $f(F_x)\cap f(F_y)\nes$, i.e. $\pl(F_x)\eta\pl(F_y)$, and, hence, by (InSkeLC), $F_x\cap F_y\nes$, a contradiction. Thus, $f(x)\neq f(y)$.
\sqs

 Again, it is easy to see that condition (InSkeLC) implies condition (InHLC) when $\p$ is a $\SKAL$-morphism.

\begin{theorem}\label{ilcompo}
Let $f\in\OLC(X,Y)$,  $(A,\rho,\BBBB)=\LAM^t(X)(=\Psi^t(X))$, $(B,\eta,\BBBB\ap)=\LAM^t(Y)(=\Psi^t(Y))$ and $\p=\LAM^t(f)(=\Psi^t(f))$. Then $f$ is an injection if and only if $\p$ is a surjection.
\end{theorem}

\doc Let $\p$ be a surjection. Let
$a,b\in \BBBB$ and $\pl(a)\eta\pl(b)$. Then, by condition (LO)
 from \cite{Di2}, $\p(\pl(a))\rho b$. Since surjectivity of $\p$ implies that $\p(\pl(a))=a$ (see
\cite[1.21]{Di2}), we get that $a\rho b$. Therefore, $\p$ satisfies condition
(InSkeLC). Hence, by Theorem \ref{ilcompnk}, $f$ is an injection.

Let $f:X\lra Y$ be an injective open map. Then
$\p(G)=f\inv(G)$, for every $G\in RC(Y)$ (see  the proof of
\cite[Theorem 2.17]{Di2}). For every $F\in RC(X)$ we have, by
\cite[Corollary 2.5]{Di2} and \cite[Lemma 2.6]{Di2}, that
$\cl(f(F))\in RC(Y)$. Set $G=\cl(f(F))$. Then, by
\cite[1.4.C]{E2}, $f\inv(G)=\cl(f\inv(f(F)))$ (because $f$ is an
open map), and the injectivity of $f$ implies that $f\inv(G)=F$.
Hence $\p(G)=F$. Therefore, $\p$ is a surjection.
\sqs

Now we will be occupied with the homeomorphic embeddings. We will  call them shortly {\em embeddings}.

 We will need a lemma from \cite{CNG}:

\begin{lm}\label{isombool}
Let $X$ be a dense subspace of a topological space $Y$. Then the
functions $r:RC(Y)\lra RC(X)$, $F\mapsto F\cap X$, and $e:RC(X)\lra
RC(Y)$, $G\mapsto \cl_Y(G)$, are Boolean isomorphisms between Boolean
algebras $RC(X)$ and $RC(Y)$, and $e\circ r=id_{RC(Y)}$, $r\circ
e=id_{RC(X)}$.
\end{lm}

\begin{theorem}\label{f1}
Let $f\in\HLC(X,Y)$, $(A,\rho,\BBBB)=\LAM^t(X)$, $(B,\eta,\BBBB\ap)=\LAM^t(Y)$ and $\p=\LAM^t(f)$. Then $f$ is a dense
 embedding iff $\p$  is a Boolean isomorphism satisfying the following condition:

\smallskip

\noindent(LO')  $\fa b\in B$ and $\fa a\in\BBBB$, $\p\inv(a)\eta b$ implies $a\rho\p(b)$.
\end{theorem}

\doc Let $f$ be a dense   embedding of $X$ in $Y$. Then
$f(X)$ is a locally compact dense subspace of $Y$ and hence it is
open in $Y$. Thus $f$ is an open injection. Therefore, by \cite[Theorem 2.17]{Di2} and \ref{mainfed}, $\p$ is a
complete homomorphism satisfying condition (LO) from \cite[Definition 2.16]{Di2}. Put $Z=f(X)$ and let $i:Z\lra Y$ be the
embedding of $Z$ in $Y$. Then $\psi=\LAM^t(i):RC(Y)\lra RC(Z)$ is
defined by the formula $\psi(F)=cl_Z(Z\cap\int_Y(F))=F\cap Z$, for
every $F\in RC(Y)$. Hence, by Lemma \ref{isombool}, $\psi$ is a
Boolean isomorphism. Since $f=i\circ f_{\upharpoonright X}$, we
obtain that $\p$ is a Boolean isomorphism as well. Then $\pl=\p\inv$ and thus condition (LO) coincides with condition (LO') (because the only one difference between the two conditions is that $\p\inv$ in (LO') is replaced with $\pl$ in (LO)).
So, $\p$ is a Boolean isomorphism satisfying condition (LO').

Conversely, let $\p$  be a Boolean isomorphism satisfying condition (LO'). Then $\p$ is a
complete homomorphism satisfying condition (LO). Obviously, condition (DLC3S) (written here immediately after \ref{dvfi}) implies condition (L1) from \cite[Definition 2.10]{Di2}.
If $a\in\BBBB$ then, by condition (DLC4) (see \ref{dvfi}), there exists $b\in\BBBB\ap$ such that $a\le\p(b)$. Then $\pl(a)\le b$. Thus $\pl(a)\in\BBBB\ap$. Therefore, $\p$ satisfies also condition (L2) from \cite[Definition 2.10]{Di2}. Hence, $\p$ is a $\OAL$-morphism (see \cite[Definition 2.16]{Di2}). Thus, by \cite[Theorem 2.17]{Di2}, $f$ is an open map. Since, by Theorem \ref{ilcompo}, $f$ is an injection, we get that $f$ is an   embedding. Finally,
by Proposition \ref{injdl},
 $f(X)$ is dense in $Y$.
\sqs

\begin{rem}\label{clemb}
\rm
Note that, in the notations of Theorem \ref{f1}, $f$ is a closed embedding iff $\p$ satisfies conditions (PAL5) (see \ref{defplc}) and (InHLC) (see \ref{injhlc}); this follows from Theorems \ref{lcper} and \ref{injhlc}.
\end{rem}

\begin{pro}\label{genemb}
Let $f\in\HLC(X,Y)$, $(A,\rho,\BBBB)=\LAM^t(X)$, $(B,\eta,\BBBB\ap)=\LAM^t(Y)$ and $\p=\LAM^t(f)$. Then $f$ is an
embedding iff there exists a complete LCA $(A_1,\rho_1,\BBBB_1)$ and $\DLC$-morphisms $\p_1:(A_1,\rho_1,\BBBB_1)\lra(A,\rho,\BBBB)$ and $\p_2:(B,\eta,\BBBB\ap)\lra(A_1,\rho_1,\BBBB_1)$ such that  $\p=\p_1\circ\p_2$, $\p_1$    is a Boolean isomorphism satisfying condition (LO')
 and $\p_2$ satisfies conditions (PAL5)  and (InHLC).
\end{pro}

\doc Obviously, $f$ is an   embedding iff $f=i\circ f_1$ where $f_1$ is a dense embedding and $i$ is a closed embedding. (Indeed, when $f$ is an embedding then let $f_1:X\lra  \cl_Y (f(X))$ be the restriction of $f$ and $i:\cl_Y(f(X))\lra Y$ be the inclusion map; the converse is also clear.) Setting $\p_1=\LAM^t(f_1)$ and $\p_2=\LAM^t(i)$, we get, by Theorem \ref{lccont}, that $\p=\p_1\circ\p_2$. Now our assertion follows from \ref{f1} and \ref{clemb}. \sqs

\begin{theorem}\label{dembz}
Let $f\in\ZLC(X,Y)$, $\p=\TE^t(f)$, $(A,I)=\TE^t(X)$ and $(B,J)=\TE^t(Y)$. Then $f$ is a dense embedding iff $\p$ is an injection and $\p(J)\spe I$.
\end{theorem}

\doc Let $f$ be a dense embedding. Then $f(X)$ is open in $Y$ and thus $f$ is an open injection. Now, Theorem \ref{injopz} implies that $\p(J)\spe I$. Since $\cl(f(X))=Y$, we get, by \ref{injdz}, that $\p$ is an injection.

Conversely, let $\p$ be an injection and $\p(J)\spe I$. Then, by \ref{injdz}, $\cl(f(X))=Y$. We will show now that $\p$ has a lower P-preadjoint. Indeed, for every $a\in I$ there exists a unique $b_a\in J$ such that $\p(b_a)=a$. Let $\psi:I\lra J$ be defined by $\psi(a)=b_a$ for every $a\in I$. Then, obviously, $\p(\psi(a))=a$, for every $a\in I$. Thus condition (OZL2) (see \ref{thzop}) is satisfied. Further, let $a\in I$, $b\in J$ and $a\we\p(b)=0$. Then $a=\p(\psi(a))$. Hence $\p(\psi(a)\we b)=0$. Since $\p$ is an injection, we get that $\psi(a)\we b=0$. So, condition (OZL1) (see \ref{thzop}) is also satisfied. Therefore, $\psi$ is a lower P-preadjoint of $\p$. Hence, by Theorem \ref{thzop}, $f$ is an open map. Now, using the condition $\p(J)\spe I$, we get, by Theorem \ref{injopz}, that $f$ is an injection. Hence, $f$ is a dense embedding. \sqs

\begin{theorem}\label{clembz}{\rm (M. Stone \cite{ST})}
Let $f\in\ZLC(X,Y)$, $\p=\TE^t(f)$, $(A,I)=\TE^t(X)$ and $(B,J)=\TE^t(Y)$. Then $f$ is a closed embedding iff $\p(J)= I$.
\end{theorem}

\doc Let $f$ be a closed embedding. Then $f$ is a perfect injection. Hence, by Theorem \ref{injperz}, $\p(J)=I$.

Conversely, let $\p(J)=I$. Then, by Theorem \ref{genstonep}, $f$ is a perfect map. Using Proposition \ref{surdz}, we get that $f$ is an injection. \sqs

\begin{pro}\label{embz}
Let $f\in\ZLC(X,Y)$, $\p=\TE^t(f)$, $(A,I)=\TE^t(X)$ and $(B,J)=\TE^t(Y)$. Then $f$ is an embedding iff there exists a ZLBA $(A_1,I_1)$ and two $\ZLBA$-morphisms $\p_1:(A_1,I_1)\lra (A,I)$ and $\p_2:(B,J)\lra (A_1,I_1)$ such that $\p=\p_1\circ\p_2$, $\p_1$ is an injection, $\p_1(I_1)\spe I$ and    $\p_2(J)= I$.
\end{pro}

\doc Arguing as in the proof  of Proposition \ref{genemb} and using Theorems \ref{dembz} and \ref{clembz}, we get the desired result. \sqs

Now we will characterize LCA-embeddings for $\SKAL$-morphisms.
This will imply  a  generalization of a theorem of Fedorchuk
\cite[Theorem 6]{F2}.

Recall that a continuous mapping $f:X\lra Y$ is said to be {\em
semi-open}\/ (\cite{V2}) if for every point $y\in f(X)$ there
exists a point $x\in f\inv(y)$ such that, for every $M\sbe X$,
$x\in\int_X(M)$ implies that $y\in\int_{f(X)}(f(M))$.

\begin{theorem}\label{lcaemb}
Let $f\in\SKLC(Y,X)$. Then $\p=\LAM^t(f)(=\Psi^t(f))$ is an
LCA-embedd\-ing iff $f$ is  a
semi-open perfect surjection.
\end{theorem}

\doc Note that, by Theorem \cite[Theorem 2.11]{Di2}, $\p$ is a complete Boolean homomorphism. Recall also that when our map $f$ is  closed then it is quasi-open (see \cite[Corollary 2.5(b)]{Di2}).

Let   $(A,\rho,\BBBB)=\LAM^t(X)$ and $(B,\eta,\BBBB\ap)=\LAM^t(Y)$.
Then  $\p:(A,\rho,\BBBB)\lra (B,\eta,\BBBB\ap)$.
Set
 $C=C_\rho$ and $C\ap=C_\eta$  (see \cite[1.16]{Di2} for
the notations). Then, by \cite[1.16]{Di2}, $(A,C)$ and $(B,C\ap)$
are CNCA's. By the proof of Theorem \cite[2.1]{Di2},
$\Psi^a(A,C)=\a X=X\cup\{\s^A_\infty\}$ and $\Psi^a(B,C\ap)=\a
Y=Y\cup\{\s^B_\infty\}$.

 Let $\p$ be an LCA-embedding, i.e. $\p:A\lra B$ is a Boolean
embedding such that, for any $a,b\in A$, $a\rho b$ iff
$\p(a)\eta\p(b)$, and $a\in\BBBB$ iff $\p(a)\in\BBBB\ap$; hence
$\p$ satisfies condition (PAL5) (see \ref{defplc}). Then, by Theorems
\ref{surplc} and \ref{lcper}, $f$ is a perfect  surjection. It remains
to show that $f$ is semi-open.
 Denote by $\p_c$ the map $\p$ regarded as a function from $(A,C)$ to
 $(B,C\ap)$.
  By \cite[(55)]{Di2}, $\p_c$ satisfies
condition (F1)  from \cite[2.12]{Di2}. We will show that $\p_c$ is
an NCA-embedding. Indeed, for any $a,b\in A$, we have that $aCb$
iff $a\rho b$ or $a,b\nin\BBBB$; since $\p$ is an LCA-embedding,
we obtain that $aCb$ iff $\p_c(a)C\ap\p_c(b)$. So, $\p_c$ is an
NCA-embedding and a $\SAC$-morphism (see \cite[Definition 2.12]{Di2}). Then, by Theorem 6 of
Fedorchuk's paper \cite{F2}, $f_c=\Psi^a(\p_c):\a Y\lra\a X$ is a
semi-open map. If $1_A\nin\BBBB$ and $1_B\nin\BBBB\ap$ then
$f_c\inv(\s^A_\infty)=\{\s^B_\infty\}$ (see the proof of
\cite[Theorem 2.15]{Di2}) and since $f=(f_c)_{|Y}$, we obtain that $f$
is semi-open.  Further, if $1_A\in\BBBB$  and $1_B\in\BBBB\ap$
then the fact that $f$ is semi-open is obvious. Since only these
two cases are possible in the given situation, we have proved that
$f$ is a perfect quasi-open semi-open surjection.

Conversely, let $f$ be a perfect semi-open surjection.
Then, by Theorem \ref{surplc}, $\p$ is an injection.
Hence $\pl\circ\p=id_A$. Thus, if
$\p(a)\in\BBBB\ap$ then, by (L2) (see \cite[Definition 2.10]{Di2}) (which follows from (DLC4)), $a=\pl(\p(a))\in\BBBB$.
Since $f$ is perfect, Theorem \ref{lcper} implies that
$\p$ satisfies condition (PLC5).
Using
it, we obtain that $a\in\BBBB$ iff $\p(a)\in\BBBB\ap$. Since
(L1) takes place (see \cite[2.10 and 2.11]{Di2}), it remains only to prove that $a\rho b$ implies
$\p(a)\eta\p(b)$, for all $a,b\in A$.  Let
$F,G\in RC(X)$,  $F\cap G\nes$ and $x\in F\cap G$. Set $U=\int(F)$
and $V=\int(G)$. Then $x\in \cl(U)\cap\cl(V)$. Since $f$ is a
semi-open surjection, there exists $y\in f\inv(x)$ such that, for
every $M\sbe Y$, $y\in\int_Y(M)$ implies that $x\in\int_X(f(M))$.
We will show that $y\in \cl(f\inv(U))\cap\cl(f\inv(V))$. Indeed,
suppose that $y\nin \cl(f\inv(U))$. Then there exists an open
neighborhood $Oy$ of $y$ such that $Oy\cap f\inv(U)=\ems$. Thus
$f(Oy)\cap U=\ems$. Since $x\in \cl(U)$ and $x\in \int(f(Oy))$, we
obtain a contradiction. Hence $y\in \cl(f\inv(U))$. Analogously we
 show that $y\in \cl(f\inv(V))$. Therefore, $y\in
\cl(f\inv(U))\cap\cl(f\inv(V))=\p(F)\cap\p(G)$. So, we get that $a\rho b$ implies
$\p(a)\eta\p(b)$. Therefore, $\p$ is an LCA-embedding.
\sqs

\section{Open sets. Regular closed  sets}

We will first recall some definitions and results from \cite{Di3}.

\begin{defi}\label{lideal}{\rm \cite{Di3}}
\rm Let $(A,\rho,\BBBB)$ be an LCA. An ideal $I$ of $A$ is called
a {\em $\d$-ideal} if $I\sbe \BBBB$ and for any $a\in I$ there
exists $b\in I$ such that $a\ll_\rho b$. If $I_1$ and $I_2$ are
two $\d$-ideals of $(A,\rho,\BBBB)$ then we put $I_1\le I_2$ iff
$I_1\sbe I_2$. We will denote by $(I(A,\rho,\BBBB),\le)$ the poset
of all $\d$-ideals of $(A,\rho,\BBBB)$.

Note that  $I$ is a $\d$-ideal of $(A,\rho_s,\BBBB)$ iff it is an ideal of $\BBBB$ (see \ref{extrcr} for $\rho_s$).
\end{defi}

\begin{fact}\label{dideal}{\rm \cite{Di3}}
Let $(A,\rho,\BBBB)$ be an LCA. Then, for every $a\in A$, the set
$\{b\in\BBBB\st b\ll_\rho a\}$ is a $\d$-ideal.
\end{fact}

The $\d$-ideals regarded in Fact \ref{dideal}
will be called\/ {\em A-principal $\d$-ideals} or, when $a\in\BBBB$, {\em principal $\d$-ideals}.
Note that the letter $``$A" coincides with the letter with which is denoted the Boolean algebra in the LCA $(A,\rho,\BBBB)$; thus, if we start with an LCA $(B,\rho,\BBBB)$, then  this kind of $\d$-ideals will be called
$``$B-principal $\d$-ideals" instead of $``$A-principal $\d$-ideals". In particular case, when $\rho=\rho_s$ (see \ref{extrcr} for $\rho_s$), we obtain the notion of {\em A-principal ideal} (see \ref{lideal}).

Recall that a {\em frame} is a complete lattice $L$ satisfying the
infinite distributive law $a\we\bigvee S=\bigvee\{a\we s\st s\in
S\}$, for every $a\in L$ and every $S\sbe L$.

\begin{fact}\label{frlid}{\rm \cite{Di3}}
Let $(A,\rho,\BBBB)$ be an LCA. Then the poset
$(I(A,\rho,\BBBB),\le)$ of all $\d$-ideals of $(A,\rho,\BBBB)$
(see \ref{lideal}) is a frame.
\end{fact}

\begin{theorem}\label{opensetsfr}{\rm \cite{Di3}}
Let $(A,\rho,\BBBB)$ be an LCA, $X=\Psi^a(A,\rho,\BBBB)$ and
$\OO(X)$ be the frame of all open subsets of\/ $X$. Then there
exists a frame isomorphism
$$\iota:(I(A,\rho,\BBBB),\le)\lra (\OO(X),\sbe),$$
where $(I(A,\rho,\BBBB),\le)$ is the frame of all $\d$-ideals of
$(A,\rho,\BBBB)$. The  isomorphism $\iota$ sends the set
$PI(A,\rho,\BBBB)$ of all A-principal $\d$-ideals of
$(A,\rho,\BBBB)$ onto the set of those regular open subsets of $X$
whose complements are in $\lag(A)$. In particular, if
$(A,\rho,\BBBB)$ is a CLCA, then $\iota(PI(A,\rho,\BBBB))=RO(X)$.
\end{theorem}

In analogy with Stone's terminology from \cite{ST1, ST}, a $\d$-ideal $J$ of an LCA
$(A,\rho,\BBBB)$  will be called a {\em simple $\d$-ideal}\/ if it has a complement in the frame $I(A,\rho,\BBBB)$, i.e. if $J\vee\neg J=\BBBB$ (here $\neg J$
is the pseudocomplement of $J$ in the frame $I(A,\rho,\BBBB)$).
(see \ref{lideal} for the notations); also, the regular elements of the frame $I(A,\rho,\BBBB)$ (i.e. those $J\in I(A,\rho,\BBBB)$ for which $\neg\neg J=J$)
will be called {\em normal $\d$-ideals}.

\begin{cor}\label{varopen}
Let $(A,\rho,\BBBB)$ be a CLCA, $(X,\OO)=\LAM^a(A,\rho,\BBBB)$, $J$ be a $\d$-ideal of $(A,\rho,\BBBB)$ and $U=\iota(J)$ (see Theorem \ref{opensetsfr} for $\iota$). Then:

\smallskip

\noindent(a) $U$ is a clopen set iff $J$ is a simple $\d$-ideal of $(A,\rho,\BBBB)$;

\smallskip

\noindent(b) $U$ is a regular open set $\iff$ $J$ is a normal $\d$-ideal of $(A,\rho,\BBBB)$ $\iff$ $J$ is an A-principal $\d$-ideal;

\smallskip

\noindent(c) $U$ is a compact open set iff $J$ is a principal ideal of\/ $\BBBB$.
\end{cor}

\doc Since the map $\iota$ is a frame isomorphism (see Theorem \ref{opensetsfr}), it preserves and reflects the regular elements and the elements which have a complement. Note also that the pseudocomplement  $\neg U$ of $U$ in the frame $(\OO,\sbe)$ is the set $\int(X\stm U)$.

\noindent(a) Clearly, $U$ is a clopen set iff it has a complement in  the frame $(\OO,\sbe)$ iff $J$ is a simple $\d$-ideal.

\noindent(b) Obviously, $U$ is a regular open set iff it is a  regular element of the frame $(\OO,\sbe)$.  Thus our assertion follows from the second statement in Theorem \ref{opensetsfr}.

\noindent(c) We have that $U$ is a compact open set $\iff$ $U=\lag(a)$ for some $a\in \BBBB$ such that $a\llx a$ $\iff$ $U=\lag(a)$ for some $a\in \BBBB$ such that the set $\{b\in\BBBB\st b\le a\}$ is a $\d$-ideal in $(A,\rho,\BBBB)$ $\iff$ $J=\{b\in\BBBB\st b\le a\}$ for some $a\in\BBBB$.
\sqs

The well-known Stone's result \cite{ST} that open sets correspond to the ideals is contained in the next corollary.

\begin{cor}\label{opensetsfrz}
Let $(A,I)$ be a ZLBA, $X=\TE^a(A,I)(=\TE^a_g(I))$ and
$(\OO(X),\sbe)$ be the frame of all open subsets of\/ $X$. Then there
exists a frame isomorphism
$$\iota:(Idl(I),\le)\lra (\OO(X),\sbe).$$
 The  isomorphism $\iota$ sends the set
$PI(A,I)$ of all A-principal ideals of
$(A,I)$ onto the set of those regular open subsets of $X$
whose complements are in $\lag(A)$. In particular, if
$A$ is a complete Boolean algebra, then $\iota(PI(A,I))=RO(X)$.
\end{cor}

\doc By  Lemma \ref{rhoslemma}, $(A,\rho_s,I)$ is an LCA. As it is noted in \ref{lideal}, in this situation the notions of $``\d$-ideal of $(A,\rho_s,I)$"  and $``$ideal of $I$" coincide. Using also the remark in \ref{dideal} about A-principal ideals, we obtain that our assertion is a particular case of Theorem \ref{opensetsfr}. Thus the isomorphism $\iota$ from Theorem \ref{opensetsfr} is the required one. \sqs

\begin{cor}\label{zvaropen}{\rm (M. Stone \cite[Theorem 5]{ST})}
Let $(A,I)$ be a ZLBA, $(X,\OO)=\TE^a(A,I)$ ($=\TE^a_g(I)$), $J$ be an ideal of $I$ and $U=\iota(J)$ (see Corollary \ref{opensetsfrz} or Theorem \ref{opensetsfr} for $\iota$). Then:

\smallskip

\noindent(a) $U$ is a clopen set $\iff$ $J$ is a simple  ideal of $I$ $\iff$ $J$ is an A-principal ideal;

\smallskip

\noindent(b) $U$ is a regular open set iff $J$ is a normal ideal of $I$;

\smallskip

\noindent(c) $U$ is a compact open set iff $J$ is a principal ideal of $I$.
\end{cor}

\doc The proofs of the statements (b) and (c), as well as of the first part of (a), are analogous to the proofs of the corresponding assertions in Corollary \ref{varopen}. So, we need only to prove the second assertion in (a).  By Proposition \ref{cox}, we have that $\lag(A)=CO(X)$. Hence, according to Corollary \ref{opensetsfrz}, A-principal ideals correspond to those regular open sets whose complements are in $\lag(A)$, i.e. in $CO(X)$. Thus, $J$ is an A-principal ideal iff $U$ is a clopen set. This assertion has also an easy direct proof. \sqs

We have seen that the open sets correspond to the $\d$-ideals. Now we are going
to describe explicitly the dual objects of the open subsets of a locally compact Hausdorff
space $X$ using only the dual object of $X$ and the corresponding $\d$-ideal.

Recall that if $A$ is a Boolean algebra and $a\in A$ then the set
$\downarrow a$ endowed with the same meets and joins as in $A$ and
with complements $b\ap$ defined by the formula $b\ap=b^*\we
a$, for every $b\le a$, is a Boolean algebra; it is denoted by
$A|a$. If $J=\downarrow (a^*)$ then $A|a$ is isomorphic to the
factor algebra $A/J$; the isomorphism $h:A|a\lra A/J$ is the
following: $h(b)=[b]$, for every $b\le a$ (see, e.g., \cite{Si2}).

\begin{theorem}\label{conopen}
Let $(A,\rho,\BBBB)$ be a CLCA, $X=\Psi^a(A,\rho,\BBBB)$ and $I$
be a $\d$-ideal of $(A,\rho,\BBBB)$. Set $a_I=\bigvee I$ and
$B=A|a_I$. For every $a,b\in B$, set $a\eta b$ iff there exist
$c,d\in I$  such that $c\rho d$, $c\le a$ and $d\le b$. Then $(B,\eta,I)$
is a CLCA. If $\p:A\lra B$ is the natural epimorphism (i.e.
$\p(a)=a\we a_I$, for every $a\in A$) then $\p$ is a $\OAL$-morphism. Set $U=\LAM^a(B,\eta,I)$ and
$f=\Psi^a(\p)$. Then $f:U\lra X$ is an open injection and
$f(U)=\iota(I)$ (see \ref{opensetsfr} for $\iota$). Hence, $U$ is
homeomorphic to $\iota(I)$ and, thus, $(B,\eta,I)$ is isomorphic to the dual object $\Psi^t(\iota(I))$ of the open set $\iota(I)$.
\end{theorem}

\doc Obviously, $B$ is a  complete Boolean algebra and
$\p$ is a surjective complete Boolean homomorphism.

Set $V=\iota(I)$ (i.e. $V=\bigcup\{\l_A^g(b)\st b\in I\}$). Then
$V$ is open in $X$ (see \ref{opensetsfr}) and
$\cl(V)=\l_A^g(a_I)$. Since $I$ is a $\d$-ideal,
$\{\int(\l_A^g(b))\st b\in I\}$ is an open cover of $V$.

If $I=\{0\}$ then $V=\ems$, $a_I=0$, $B=\{0\}$ and $U=\ems$;
hence, in this case the assertion of the theorem is true. Thus,
let us assume that $I\neq \{0\}$.

We will first check that $(B,\eta,I)$ is a CLCA, i.e. that
conditions  (C1)-(C4)  and (BC1)-(BC3) from \cite{Di2} are
fulfilled. Note that, for every $a,b\in B$, $a\eta b$ implies that $a\rho b$; thus, if $a\llx b$ then $a(-\rho) b^*$ and hence $a(-\rho) (b^*\we a_I)$, which implies that $a\lle b$.

Let $b\in B\stm \{0\}$. Then $b=\bv\{c\we b\st c\in I\}$. Thus there exists $c\in I$ such that $c\we b\neq 0$. We get that $d=c\we b\in I$, $d\le b$ and $d\rho d$. Therefore $b\eta b$.
 So, the axiom (C1) is fulfilled. Using the same notations, we get that there exists $a\in\BBBB\stm\{0\}$ such that $a\llx d$. Then $a\in I\stm\{0\}$ and $a\lle b$.
Therefore, the axiom (BC3) is checked as well. Clearly, the axioms
(C2), (C3) and (BC2) are satisfied. Let $a,b_1,b_2\in B$ and $a\eta(b_1\vee b_2)$. Then there exist $c,d\in I$ such that $c\le a$, $d\le b_1\vee b_2$ and $c\rho d$. Since $d=(d\we b_1)\vee(d\we b_2)$, we get that either $c\rho (d\we b_1)$ or $c\rho (d\we b_2)$. Clearly, this implies that either $a\eta b_1$ or $a\eta b_2$. The converse implication is obvious. So,
we obtain that the axiom (C4) is also fulfilled.

 Let $a\in I$, $b\in B$ and $a\ll_\eta b$.  Then $a(-\eta)(b^*\we a_I)$. Thus, for every $c\in I$ such that $c\le b^*$, we have that $a(-\rho)c$. Since $a\in I$ and $I$ is a $\d$-ideal, we get that there exists $c\in I$ such that $a\llx c\llx a_I$. Then $c\we b^*\le b^*$ and $c\we b^*\in I$. Thus $a(-\rho)(c\we b^*)$, i.e $a\llx (c^*\vee b)$. Combining this fact with the inequality $a\llx c$, we get that $a\llx (c\we(c^*\vee b))$, i.e. $a\llx (b\we c)$. Then there exists $d\in\BBBB$ such that $a\llx d\llx (c\we b)$. Since $c\we b\in I$, we get that $d\in I$. Therefore, $a\llx d\llx b$. This implies that $a\lle d\lle b$.
  Thus, the
axiom (BC1) is checked.

So, we have proved that $(B,\eta,\BBBB\ap)$ is a CLCA.

We will show that $\p$ is a $\OAL$-morphism, i.e. that $\p$ satisfies axioms (L1), (L2) and (LO) from
\cite{Di2}. Note first that, for every $a\in B$,
\begin{equation}\label{eq1}
\pl(a)=a.
\end{equation}
This observation  shows that $\p$ satisfies conditions (L2) and (EL1) from \cite{Di2} (note that condition (EL1) is equivalent to
the condition (L1)). Let
us prove that the axiom (LO) is fulfilled as well. Let $a\in A$,
$b\in I$ and $\pl(b)\rho a$. Then $a\rho b$. We have to show that $b\eta\p(a)$, i.e. that $b\eta(a\we a_I)$. Suppose that
$b(-\eta)(a\we a_I)$. Then, for every $c\in I$ such that $c\le a$, we have that $b(-\rho)c$. Since $I$ is a $\d$-ideal, there exists $d\in I$ such that $b\llx d\llx a_I$. Then $d\we a\in I$ and $d\le a$. Hence $b(-\rho)(d\we a)$, i.e. $b\llx(d^*\vee a^*)$. Since $b\llx d$, we get that $b\llx(d\we(d^*\vee a^*))$. Thus $b\llx a^*$, i.e. $b(-\rho)a$, a contradiction.
Therefore,  condition (LO) is checked. So,
$\p$ is a $\OAL$-morphism. Since $\p$ is a surjection, Theorem
\ref{ilcompo} implies that $f:U\lra X$ is an open injection and hence $f$ is a
homeomorphism between $U$ and $f(U)$. Let us show that $f(U)=V$.
Recall that the function $f$ is defined by the formula $f(\s_u)=\s_{\p\inv(u)}(=\s_{\pl(u)})$, where $u$ is a bounded ultrafilter in $(B,\eta,I)$. We have also that $V=\bigcup\{\lag(a)\st a\in I\}$. Now, if $\s_u\in U$ then there exists $a\in I\cap u$. Since $\p(a)=a\we a_I=a$, we get that $a\in\p\inv(u)$. Thus $f(\s_u)\in\lag(a)\sbe V$. Hence, $f(U)\sbe V$. Conversely, if $\s\ap\in V$ then 
there exists $a\in\s\ap\cap I$. Thus, there exists an ultrafilter $v$ in $A$ such that $a\in v$ and $\s\ap=\s_v$. 
Obviously, $v\cap I$ is a filter-base of $v$ (because $a\in v\cap I$ and $I$ is an ideal). It is clear that $u=v\cap B$ is a bounded filter in $(B,\eta,I)$. Moreover, $u$ is an ultrafilter in $B$. Indeed, let $c\in B=\downarrow_A(a_I)$. If $c\in v$ then $c\in u$. If $c^*\in v$ then $a\we c^*\in v\cap B$ and thus $c^*\we a_I\in u$, i.e. $c\ap\in u$. Hence, $u$ is a bounded ultrafilter in $(B,\eta,I)$. Since $\pl(u)=u$ and $u$ is a filter-base of $v$, we get that $f(\s_u)=\s_{\pl(u)}=\s_v=\s\ap$. Therefore, $f(U)=V$. \sqs

The analogous question for the dualities between the categories $\ZLC$ and $\ZLBA$, and between $\PZLC$ and $\GBPL$, is much easier. Since the category $\PZLC$ is a cofull subcategory of the category $\ZLC$, it is enough to describe the dual $\GBPL$-objects of the objects of the category $\ZLC$, as it is done by M. Stone in \cite{ST}. If $I$ is the dual object of some $X\in\card{\ZLC}$ then $(Si(I),e_I(I))$ is its dual object in $\ZLBA$. Since for every $X\in\card{\ZLC}$  and every open subset $U$ of $X$ we have that $CK(U)=\{F\in CK(X)\st F\sbe U\}=I_U=\iota\inv(U)$ and $I_U$ is an ideal of $CK(X)$, we get that $I_U$ is not only the ideal of $CK(X)$ corresponding to $U$ but it is also the dual object of $U$ in the category $\GBPL$, i.e. $I_U=\TE^t_g(U)$; thus the dual object $\TE^t(U)$ of $U$ in $\ZLBA$ is $(Si(I_U),e_{I_U}(I_U))$. Conversely, if $I$ is an ideal of $CK(X)$ then $I=\iota\inv(\iota(I))$ and hence $I=\TE^a_g(\iota(I))$.

We will now show how one can build the CLCAs corresponding to the
regular closed subsets of a locally compact Hausdorff space $Y$
from the CLCA $\Psi^t(Y)$.

\begin{theorem}\label{conregclo}
Let $(A,\rho,\BBBB)$ be a CLCA, $X=\Psi^a(A,\rho,\BBBB)$, $a_0\in
A$ and $F=\l_A^g(a_0)$. Set $B=A|a_0$ and let $\p:A\lra B$ be the
natural epimorphism (i.e. $\p(a)=a\we a_0$, for every $a\in A$).
Put $\BBBB\ap=\p(\BBBB)$ and let, for every $a,b\in B$, $a\eta b$
iff $a \rho b$. Then $(B,\eta,\BBBB\ap)$ is a CLCA. If
$G=\Psi^a(B,\eta,\BBBB\ap)$ then $G$ is
homeomorphic to $F$ and thus $\Psi^t(F)$ is isomorphic to $(B,\eta,\BBBB\ap)$. If $f=\Psi^a(\p)$ then $f:G\lra X$ is
a closed quasi-open injection and $f(G)=F$.  
\end{theorem}

\doc We have that $B$
is a complete Boolean algebra, $\p$ is a complete Boolean
homomorphism and $\pl(a)=a$, for every $a\in B$. Set
$\psi=\l_A^g\circ\pl$. We will show that
$\psi\ap=\psi_{\upharpoonright B}$ is a Boolean isomorphism between
$B$ and $RC(F)$. Since $F\in
RC(X)$, we have, as it is well known, that $RC(F)\sbe RC(X)$ and
$RC(F)=\{G\we F\st G\in RC(X)\}$; moreover, $RC(F)=RC(X)|F$.
 Hence $\psi\ap:B\lra RC(F)$ is a Boolean isomorphism.  For any
$a,b\in B$, we have that $a\eta b\iff \pl(a)\rho\pl(b)\iff
\l_A^g(\pl(a))\cap\l_A^g(\pl(b))\nes\iff \psi(a)\rho_F\psi(b)$.
Finally,  for any $a\in B$, we have that $a\in\BBBB\ap\iff
\pl(a)\in\BBBB\iff \l_A^g(\pl(a))$ is compact $\iff \psi(a)\in
CR(F)$. Therefore, $(B,\eta,\BBBB\ap)$ is a CLCA because $(RC(F),\rho_F,CR(F))$ is such, and they are isomorphic.
 For showing that $f:G\lra X$
is a homeomorphic embedding and $f(G)=F$, note that $\p$ satisfies
conditions (L1)-(L3) from \cite{Di2} and condition (InSkeLC), and
hence, by Theorem \ref{ilcompnk}, $f$ is a quasi-open perfect
injection, i.e. $f$ is a homeomorphic embedding. From
\cite[(45)]{Di2} we get that, for every $b\in\BBBB\ap$,
$f(\l_B^g(b))=\l_A^g(\pl(b))=\l_A^g(b)\sbe F$. Since
$G=\bigcup\{\l_B^g(b)\st b\in\BBBB\ap\}$, we obtain that $f(G)\sbe
F$. Let $y\in\int_X(F)$. Then there exists $b\in\BBBB$ such that
$y\in\int(\l_A^g(b))\sbe \l_A^g(b)\sbe\int_X(F)$. Hence
$b\in\BBBB\ap$. Using again \cite[(45)]{Di2}, we get that $y\in
f(\l_B^g(b))$, i.e. $y\in f(G)$. Thus $\int_X(F)\sbe f(G)$. Since
$f(G)$ is closed in $X$, we conclude that $f(G)\spe F$. Therefore,
$f(G)=F$.
\sqs

Now, for every $X\in\card{\ZLC}$, we will describe the dual objects $\TE^t(F)$ and $\TE^t_g(F)$ of the closed or regular closed subsets of $X$ by means of the dual objects $\TE^t(X)$ and $\TE^t_g(X)$ of $X$. The obtained result for regular closed subsets of $X$ seems to be new even in the compact case.

\begin{theorem}\label{clregclz}
Let $I,J\in\card{\GBPL}$, $X=\TE^a_g(I)$ and $F=\TE^a_g(J)$. Then:

\smallskip

\noindent(a){\rm (M. Stone \cite[Theorem 4(4)]{ST})} $F$ is homeomorphic to a closed subset of $X$ iff there exists a  0-pseudolattice epimorphism $\p:I\lra J$ (i.e. iff $J$ is a quotient  of $I$);

\smallskip

\noindent(b) $F$ is homeomorphic to a regular closed subset of $X$ if and only if there exists a  0-pseudolattice epimorphism $\p:I\lra J$ which preserves all meets that happen to exist in $I$.
\end{theorem}

\doc (a) Let $F$ be homeomorphic to a closed subset of $X$, i.e there exists a closed embedding $f:F\lra X$. Then, by Theorem \ref{clembz}, $\p\ap=\TE^t_g(f):\TE^t_g(X)\lra\TE^t_g(F)$ is a surjective 0-pseudolattice homomorphism. Thus, by the duality, there exists a a surjective 0-pseudolattice homomorphism
$\p:I\lra J$.

Conversely, if $\p:I\lra J$ is a surjective 0-pseudolattice homomorphism then, by Theorem \ref{clembz},  $F$ is homeomorphic to a closed subset of $X$.

\smallskip

\noindent(b) Having in mind the assertion (a) above and Theorem \ref{thqgbpl}, it is enough to show that if $f:F\lra X$ is a closed injection then $f(F)\in RC(X)$ iff $f$ is a quasi-open map. This can be easily proved, so that the proof of assertion (b) is complete. \sqs

We will finish with mentioning some assertions about isolated points.
All these statements have  easy proofs which will be omitted.

\begin{pro}\label{isopo}
(a) Let $(A,\rho,\BBBB)$ be an LCA, $X=\Psi^a(A,\rho,\BBBB)$ and $a\in
A$. Then $a$ is an atom of $A$ iff  $\l_A^g(a)$ is an isolated
point of the space $X$. Also, for every isolated point $x$ of $X$ there exists an
$a\in \BBBB$ such that $a$ is an atom of\/ $\BBBB$ (equivalently, of $A$) and  $\{x\}=\lag(a)$.

\smallskip

\noindent(b) Let $(A,I)$ be a ZLBA, $X=\TE^a(A,I)$ and $a\in
A$. Then $a$ is an atom of $A$ iff  $\l_A^g(a)$ is an isolated
point of the space $X$. Also, for every isolated point $x$ of $X$ there exists an
$a\in I$ such that $a$ is an atom of $I$ (equivalently, of $A$) and  $\{x\}=\lag(a)$.
\end{pro}

\begin{pro}\label{disspa}
(a) Let $(A,\rho,\BBBB)$ be an LCA and $X=\Psi^a(A,\rho,\BBBB)$. Then
$X$ is a discrete space  iff\/ $\BBBB$ coincides with the set of all finite
sums of the atoms of $A$.

\smallskip

\noindent(b) Let $(A,I)$ be a ZLBA and $X=\TE^a(A,I)(=\TE^a_g(I))$. Then
$X$ is a discrete space  $\iff$ $I$ coincides with the set of all finite
sums of the atoms of $A$ $\iff$ the elements of $I$ are either atoms of $I$ or finite sums of atoms of $I$.
\end{pro}

\begin{pro}\label{extrspa}
(a) Let $(A,\rho,\BBBB)$ be a CLCA and $X=\Psi^a(A,\rho,\BBBB)$. Then
$X$ is an extremally disconnected  space  iff $a\llx a$, for every
$a\in A$.

\smallskip

\noindent(b){\rm (M. Stone \cite{ST})} Let $(A,I)$ be a ZLBA and $X=\TE^a(A,I)(=\TE^a_g(I))$. Then
$X$ is an extremally disconnected  space  iff $A$ is a complete Boolean algebra.
\end{pro}

\begin{pro}\label{denseiso}
(a) Let $(A,\rho,\BBBB)$ be an LCA and $X=\Psi^a(A,\rho,\BBBB)$. Then
the set of all isolated points of $X$ is dense in $X$ iff $A$ is
an atomic Boolean algebra iff $\BBBB$ is an atomic 0-pseudolattice.

\smallskip

\noindent(b) Let $(A,I)$ be a ZLBA and $X=\TE^a(A,I)(=\TE^a_g(I))$. Then
the set of all isolated points of $X$ is dense in $X$ iff $A$ is
an atomic Boolean algebra iff $I$ is an atomic 0-pseudolattice.
\end{pro}

\baselineskip = 0.75\normalbaselineskip

\end{document}